\documentclass{elsarticle}

\usepackage{lineno,hyperref}
\usepackage{a4wide}
\usepackage{amsmath}
\usepackage{amsthm}
\usepackage{amssymb}
\usepackage{amsfonts}
\usepackage{mathrsfs}
\usepackage{stmaryrd}
\usepackage{booktabs}
\usepackage{wrapfig}
\usepackage{graphicx}
\usepackage{epstopdf}
\usepackage{paralist}
\usepackage{subcaption}
\usepackage{xcolor}
\usepackage{enumitem}
\usepackage{siunitx}
\usepackage{multirow}
\usepackage{xspace}

\usepackage[capitalise]{cleveref}
\usepackage{todonotes}
\usepackage{bm}

\newcommand{\R}{\mathbb{R}}

\newcommand{\calD}{\mathcal{D}}
\newcommand{\calE}{\mathcal{E}}

\newcommand{\calP}{\mathcal{P}}

\newcommand{\calT}{\mathcal{T}}
\newcommand{\calV}{\mathcal{V}}
\newcommand{\calS}{\mathcal{S}}
\newcommand{\calM}{\mathcal{M}}
\newcommand{\upFactor}{\zeta}

\newcommand{\n}{\mathbf{n}}

\newcommand{\meas}[1]{|{#1}|}

\newcommand{\ff}{\mathrm{ff}}
\newcommand{\porm}{\mathrm{pm}}
\newcommand{\ipmff}{\mathrm{if}}
\newcommand{\vel}{\mathbf{v}}
\newcommand{\BFCoeff}{\alpha_{\mathrm{BJS}}}
\newcommand{\vertex}{v}

\newcommand{\up}{\mathrm{up}}
\newcommand{\dn}{\mathrm{dn}}

\newcommand{\error}{e}

\newcommand{\eg}{e.g. }

\newcommand{\ie}{i.e.\ }

\newcommand{\BOX}{\textsc{box}\xspace}

\newcommand{\elem}{E}
\newcommand{\gridFace}{\epsilon}
\newcommand{\scv}{\varkappa}
\newcommand{\discFlux}{F}
\newcommand{\basisFunc}{\Phi}
\newcommand{\iSection}{\gamma}
\newcommand{\iSectionSet}{\mathcal{I}}

\newcommand{\someQuantity}{s}
\newcommand{\someFlux}{\boldsymbol{\varphi}}

\newcommand{\adv}{\mathrm{adv}}
\newcommand{\diff}{\mathrm{diff}}
\newcommand{\cond}{\mathrm{cond}}

\theoremstyle{plain}

\newtheorem{definition}{Definition}

\theoremstyle{remark}

\newcommand{\Dumux}{{Du\-Mu$^\text{x}$ }}
\newcommand{\eqbydef}{:=}
\newcommand{\ltwoproj}{\ensuremath{\Pi^{\mathrm{L}_2}}\xspace}
\newcommand{\avgproj}{\ensuremath{\Pi^{\Gamma}}\xspace}

\makeatletter
\def\ps@pprintTitle{%
 \let\@oddhead\@empty
 \let\@evenhead\@empty
 \def\@oddfoot{}%
 \let\@evenfoot\@oddfoot}
\makeatother

\begin{document}

\begin{frontmatter}

\title{Coupling staggered-grid and vertex-centered finite-volume methods for coupled porous-medium free-flow problems}


\author[addressIWS]{Martin Schneider\corref{mycorrespondingauthor}}
\ead{martin.schneider@iws.uni-stuttgart.de}
\author[addressIWS]{Dennis Gl\"aser}
\ead{dennis.glaeser@iws.uni-stuttgart.de}
\author[addressIWS]{Kilian Weishaupt}
\ead{kilian.weishaupt@iws.uni-stuttgart.de}
\author[addressIWS]{Edward Coltman}
\ead{edward.coltman@iws.uni-stuttgart.de}
\author[addressIWS]{Bernd Flemisch}
\ead{bernd.flemisch@iws.uni-stuttgart.de}
\author[addressIWS]{Rainer Helmig}
\ead{rainer.helmig@iws.uni-stuttgart.de}

\cortext[mycorrespondingauthor]{Corresponding author}

\address[addressIWS]{Institute for Modelling Hydraulic and Environmental Systems,
           University of Stuttgart,
             Pfaffenwaldring 61,
               70569 Stuttgart, Germany}

\begin{abstract}
In this work, a new discretization approach for coupled free and porous-medium flows is introduced,
which uses a finite volume staggered-grid method for the discretization of the Navier-Stokes equations
in the free-flow subdomain, while a vertex-centered finite volume method is used in the porous-medium
flow domain. The latter allows for the use of unstructured grids in the porous-medium subdomain,
and the presented method is capable of handling non-matching grids at the interface.
In addition, the accurate evaluation of coupling terms and of additional nonlinear velocity-dependent terms in the porous medium is ensured by the use of basis functions and by having degrees of freedom naturally located at the interface.
The available advantages of this coupling method are investigated in a series of tests: a convergence test for various grid types, an evaluation of the implementation of coupling conditions, and an example using the velocity dependent Forchheimer term in the porous-medium subdomain.
\end{abstract}


\end{frontmatter}


\section{Introduction}
\label{sec:introduction}
Descriptions of flow and transport processes in our natural, biological, and engineered environment exist for a diverse spectrum of flow conditions.
Isolated from each other, these descriptive models can be used to effectively simulate flow within their domain of analysis.
On the contrary, most relevant applications contain multiple coupled flow regimes.
One common example of multi-domain systems that are challenging to resolve with an isolated flow model are coupled free-flow and multi-phase porous-medium flows.
Mass, momentum, and heat exchange will occur across the interface between these domains dependent on the conditions in each subdomain.
The flow fields on either side of the interface will change with this exchanged content, further driving the exchange processes between them.

Applications that rely on a detailed description of this exchange are abundant, ranging from large scale environmental applications such as soil-water evaporation and salt precipitation \cite{vanderborght2017a,jambhekar2015a}, to the small scale exchange processes within a PEM fuel cell \cite{gurau2009a}.
Evaluating these applications via direct numerical simulation of the flow fields, like those shown in \cite{chu2021a}, is usually infeasible due to the immense computational cost and the effort involved in procuring pore geometries.
For the analysis of most applications, evaluations using the concept of a representative elementary volume (REV) scale \cite{bear1972a} can be used, where averaging techniques create upscaled descriptions of flow in a porous domain \cite{whitaker1999a}.

Evaluating coupled free-flow and porous-medium flow systems at this REV scale can then be done using a single- or two-domain approach.
A single-domain approach evaluates both domains with a single set of equations, and is possible using the Brinkman equations \cite{brinkman1949a}; this method has been expanded upon significantly in \cite{Carrillo2020a}. With a two-domain approach, momentum conservation in each domain can be described using a model designed for its context, in this case, the Navier-Stokes equations for describing free flow, and the multi-phase Darcy's law for describing multi-phase flow in a porous medium \cite{mosthaf2011a}.
In order to couple these two flow domains, thermodynamically consistent interface conditions have to be used to enforce mass, momentum, and energy exchanges between these two flow domains \cite{hassanizadeh1989a}.

Numerous discretization methods for these coupled two-domain models have been developed.
Some models, e.g. \cite{fetzer2017a, schneider2020}, couple cell-centered finite volumes, used to discretize the porous medium, with a Marker-and-Cell (staggered grid) finite volume scheme \cite{harlow1965a}, used for the discretization of the free-flow domain.
Other approaches utilize the same discretization method in each domain, using either staggered grids \cite{iliev2004a}, finite elements \cite{discacciati2009a}, or a vertex-centered finite volume method, otherwise known as the \BOX scheme \cite{mosthaf2011a, huber2000a} for both domains.
Other models have been developed to use an extra lower-dimensional grid, or interface mortar elements, at the interface where the exchange fluxes can be calculated \cite{baber2016a,boon2020a}.

For many porous media based applications, variations of the core mathematical models have been developed to describe specific physical phenomena by considering additional terms in the balance equations.
These terms include the Forchheimer term, describing viscous stresses within a porous medium \cite{nield2006a}, and dispersion terms, describing transport due to subscale fluctuations in velocity \cite{scheidegger1961a}.
Both of these terms require a reconstruction of velocity for collocated discretization methods, typically requiring an expansion of local stencils \cite{srinivasan2013}.
Vertex-centered finite volume methods, such as the \BOX scheme, show their strength in this reconstruction, where the basis functions can be used to determine the velocities used in these terms \cite{huber2000a}.
In this work, the free flow is modeled using the Navier-Stokes equations discretized with the staggered-grid method.
The porous-medium flow is modeled using the multi-phase Darcy's law including the Forchheimer term, and is discretized using a vertex-centered finite volume method.

The structure of this work is as follows: \cref{sec:equations} reviews the governing model equations, \cref{sec:discretization} gives an overview of the discretization methods used in both domains, and \cref{sec:numresults} examines the performance of the proposed
coupling method by means of numerical test cases.
A review of this work, and an outlook towards further goals is provided in \cref{sec:conclusions}.

\section{Governing Equations}
\label{sec:equations}

The following section provides an outline of the assumptions made regarding geometry while developing this work, and reviews the notation.
Furthermore, the governing equations for the free-flow and the porous-medium flow domains are presented, followed by the
conditions used to couple them.

The domain in question, $\Omega\subset\R^d$, $d\in \{2, 3\}$, is an open, connected Lipschitz domain with boundary $\partial\Omega$ and $d$-dimensional measure $\meas{\Omega}$.
This domain $\Omega$ is then partitioned disjointly to $\Omega^\ff$ and $\Omega^\porm$ representing the free-flow and porous-medium flow subdomains, respectively.
The boundaries of each subdomain are defined as the interface $\Gamma^\ipmff \eqbydef \partial \Omega^\ff \cap \partial \Omega^\porm$ as well as the remainders $\Gamma^\porm \eqbydef \partial \Omega^\porm \setminus \Gamma^\ipmff$ and $\Gamma^\ff \eqbydef \partial \Omega^\ff \setminus \Gamma^\ipmff$. In the following, the superscripts are omitted if no ambiguity arises.

Using a subscript to denote a prescribed boundary condition for either pressure, $p$, or velocity, $\mathbf{v}$, the external boundary $\partial \Omega = \Gamma^\ff \cup \Gamma^\porm$ is further decomposed such that $\Gamma^\ff = \Gamma_\mathrm{v}^\ff \cup \Gamma_\mathrm{p}^\ff$ and $\Gamma^\porm = \Gamma_\mathrm{v}^\porm \cup \Gamma_\mathrm{p}^\porm$ disjointly.
Further, we assume that for a subset of the boundary $\partial \Omega$, a pressure boundary condition is imposed with positive measure, i.e. $\meas{\Gamma_\mathrm{p}^\porm \cup \Gamma_\mathrm{p}^\ff} > 0$, ensuring that the resulting system can be uniquely solved.
The vector $\n$ denotes the outwardly oriented unit normal vector on $\partial \Omega$ with respect to $\Omega$.

Molecular diffusion and heat conductions are modeled after the following laws;
\begin{align}
&\text{Fick's law: } \quad &&\mathbf{J}_{\alpha}^{\kappa, \diff} := -\varrho_{\alpha} D_{\alpha}^\kappa \nabla X_\alpha^\kappa, \label{eq:FicksLaw}\\
&\text{Fourier's law: } \quad &&\mathbf{J}^\cond := -\lambda_{\text{eff}} \nabla T, \label{eq:FouriersLaw}
\end{align}

where $\varrho_{\alpha}$ is the mass density of phase $\alpha$, $D_\alpha^\kappa$ is the diffusion coefficient of component $\kappa$ in phase $\alpha$, and $X_\alpha^\kappa$ is the mass fraction of component $\kappa$ in phase $\alpha$. $\lambda_\text{eff}$ is the effective thermal conductivity \cite{vanderborght2017a}, which in $\Omega^\ff$ is equal to the gas-phase thermal conductivity, and $T$
is the temperature, assuming local thermal equilibrium in $\Omega^\porm$.

These laws form an essential part of the governing equations in the subdomains and at the interface,
which are presented in the subsequent sections.

\subsection{Free Flow}
\label{sub:free_flow}
Within the coupled model, the free-flow subdomain consists of a gas phase with multiple components $\kappa$. Here, we drop the subscript $g$ which relates all quantities to the gas phase. The mass balance \eqref{eq:ffmass} for each component $\kappa$, the energy balance \eqref{eq:ffenergy}, and the momentum balance equations \eqref{eq:ffmomentum} are given as:
\begin{subequations}
  \label{eq:navierstokes}
  \begin{align}
    \frac{\partial \left(\varrho X^\kappa\right)}{\partial t}
   + \nabla \cdot \{ \varrho  X^\kappa \mathbf{v}
    + \mathbf{J}^{\kappa, \diff} \}
      &= q^\kappa, && \label{eq:ffmass}\\
     \frac{\partial (\varrho  u)}{\partial t}
     + \nabla \cdot \{ \varrho h \mathbf{v} + \sum_{\kappa} h^{\kappa} \mathbf{J}^{\kappa, \diff}
     + \mathbf{J}^\cond \}  &= q_T , && \label{eq:ffenergy}\\
       \frac{\partial (\varrho \vel)}{\partial t} + \nabla \cdot \{\varrho \vel \vel^{\mathrm{T}}
      - \mu (\nabla \vel + (\nabla \vel)^{\mathrm{T}})
      + p \mathbf{I} \}
      &= \varrho \textbf{g} + \mathbf{f}, &&\mathrm{in} \, \Omega^\ff, \label{eq:ffmomentum}
   \end{align}
\end{subequations}
together with appropriate boundary conditions.
The unknown variables are the gas-phase velocity $\vel$, pressure $p$ and temperature $T$, and one of the mass $X^\kappa$ or mole fractions $x^\kappa$. Here, $\varrho$ and $\mu$ denote the potentially solution-dependent gas-phase density and viscosity, respectively, while $q$ are source (or sink) terms, and $\mathbf{g}$ describes the influence of gravity. $\mathbf{f}$ is an additional momentum sink or source term. $\mathbf{I}$ is the identity tensor in $\R^{d \times d}$.
$h$ and $u$ denote the enthalpy and internal energy, respectively.

\subsection{Porous-Medium Flow}
\label{sub:porous_medium_flow}

The equations governing compositional non-isothermal two-phase (water and air) flow in the porous-medium subdomain, $\Omega^\porm$, are given by:
\begin{subequations}
    \label{eq:darcy}
  \begin{align}
 \phi \frac{\partial (\sum_\alpha \varrho_{\alpha} X_\alpha^\kappa S_\alpha )}{\partial t} + \sum\limits_\alpha \nabla \cdot  \{ \varrho_{\alpha} X_\alpha^\kappa \mathbf{v}_\alpha  + \mathbf{J}_\alpha^{\kappa, \diff} \}
  &= \sum_{\alpha} q_{\alpha}^\kappa, && \label{eq:darcyMassBalance}\\
\phi \frac{ \partial \left( \sum_\alpha \rho_\alpha u_\alpha S_\alpha \right) }{\partial t} + \left( 1 - \phi \right) \frac{\partial \left( \rho_s c_s T \right)}{\partial t} + \sum_\alpha \nabla \cdot \{ \rho_\alpha h_\alpha \mathbf{v}_\alpha \}  + \nabla \cdot \mathbf{J}^\cond
          &= q_T, && \label{eq:darcyEnergyBalance}\\
\vel_\alpha
+ \frac{k_{r \alpha}}{\mu_\alpha} \mathbf{K} \nabla \left(p_\alpha
+ \rho_\alpha \mathbf{g} \right) &=  0,  &&\mathrm{in} \, \Omega^\porm,\label{eq: darcys law}
  \end{align}
\end{subequations}
together with appropriate boundary conditions. $S_\alpha$ is the saturation of phase $\alpha$ while $\phi$ is the porosity of the permeable medium. $\varrho_s$ and $c_s$ are the solid-phase density and heat capacity.

\cref{eq: darcys law} states that the momentum balance in the porous medium is given by Darcy's law, with $\mathbf{K}$ being the permeability tensor. Similar to the free-flow equations \eqref{eq:navierstokes}, $q$ denotes source or sink terms.
The relative permeability, $k_{r \alpha}(S_w)$, as well as the capillary pressure, $p_c(S_w) = p_g - p_w$, used as closure relations, are modeled using the van Genuchten approach \cite{vanGenuchten1980,Luckner1989}, where $g$ and $w$ denote the gas and water phase, respectively.

\subsection{Coupling Conditions}
\label{sub:coupling}

In order to derive a thermodynamically consistent formulation of the coupled problem,
the conservation of mass, energy, and momentum must be guaranteed at the interface between the porous
 medium and the free-flow domain. We therefore impose the following interface conditions:
\begin{subequations}\label{eq:interfaceConditions}
 \begin{align}
 [X_g^\kappa]^\ff - [X_g^\kappa]^\porm &= 0, \label{eq:ifMassFracCont}\\
    \left[\left(\varrho_g X_g^\kappa \vel_g + \mathbf{J}^{\kappa, \diff} \right)\cdot \n \right]^\ff + \left[\sum_\alpha \left(\varrho_\alpha X_\alpha^\kappa \vel_\alpha + \mathbf{J}^{\kappa,  \diff}_{\alpha} \right)\cdot \n \right]^\porm & = 0, \label{eq:ifMassFluxCont}\\
    [T]^\ff - [T]^\porm &= 0, \label{eq:ifTempCont}\\
   \left[ \left(\varrho_g h_g \vel_g + \sum_{\kappa} h^{\kappa}_g  \mathbf{J}^{\kappa, \diff} + \mathbf{J}^\cond \right)\cdot \n \right]^\ff + \left[\left( \sum_\alpha  \varrho_\alpha h_\alpha \vel_\alpha + \mathbf{J}^\cond \right)\cdot \n \right]^\porm & = 0, \label{eq:ifTempFluxCont} \\
    \n \cdot \left[\varrho_g \vel_g \vel_g^\mathrm{T} - \mu_g (\nabla \vel_g + (\nabla \vel_g)^\mathrm{T}) + p_g \mathbf{I} \right]^\ff \n &=
      p_g^\porm, \label{eq:ifNormalMomentum} \\
     \left( - (\nabla \vel_g + (\nabla \vel_g)^\mathrm{T}) \mathbf{n} - \frac{\BFCoeff}{\sqrt{\mathbf{t} \cdot \mathbf{K} \mathbf{t}}} \vel_g \right)^\ff \cdot \mathbf{t} &=  0. \label{eq:ifBJS}
  \end{align}
\end{subequations}

The momentum transfer normal to the interface is given by \eqref{eq:ifNormalMomentum} \cite{layton2002}.
Condition \eqref{eq:ifBJS} is the commonly used Beavers-Joseph-Saffman slip condition \cite{beavers1967a, saffman1971a}. Here, $\mathbf{t}$ denotes any unit vector from the tangent bundle of $\Gamma^\ipmff$ and $\BFCoeff$ is a parameter to be determined experimentally. We remark that this condition is technically a boundary condition for the free flow, not a coupling condition between the two flow regimes. Furthermore, it has been developed for free flow strictly parallel to the interface and might lose its validity for other flow configurations \cite{yang2019a}. We are nevertheless going to use \cref{eq:ifBJS}  in the upcoming numerical examples which also include non-parallel flow at the interface,
thus accepting a potential physical inconsistency for the sake of numerical verification. \Cref{eq:ifBJS} could be replaced by a different condition in future work, such as those outlined in \cite{eggenweiler2021a,sudhakar2021a}, as this condition is not essential to the numerical model.


\section{Discretization}
\label{sec:discretization}

This section provides a brief overview of the notation regarding the partition of $\Omega$,
before the numerical schemes used in the subdomains and the incorporation
of the coupling conditions \eqref{eq:interfaceConditions} are outlined.
\begin{definition}[Grid discretization]
\label{def:griddisc}
The tuple \mbox{$\calD \eqbydef (\calM, \calS, \calT,\calE,\calP,\calV)$} denotes the grid discretization, in which
  \begin{enumerate}[label=(\roman*)]
  \item $\calM$ is the set of primary grid elements such that
    \mbox{$\overline{\Omega}= \cup_{\elem \in \calM} \overline{\elem}$}.
    For each element $\elem\in\calM$, $\meas{\elem}>0$ and $\mathbf{x}_\elem$ denote its volume and barycenter.
  \item $\calS$ is the set of faces such that each face $\gridFace$ is a $(d-1)$-dimensional hyperplane
        with measure
        \mbox{$\meas{\gridFace}>0$}.
        For each primary grid element $\elem \in \calM$, $\calS_\elem$ is the subset of
        $\calS$ such that \mbox{$\partial \elem  = \cup_{\gridFace \in\calS_\elem}{\gridFace}$}. Furthermore, $\mathbf{x}_\gridFace$ denotes
        the barycenter and $\n_{\elem,\gridFace}$ the unit vector that is normal to
        $\gridFace$ and outward to $\elem$.
  \item $\calT$ is the set of control volumes (dual grid elements) such that
    \mbox{$\overline{\Omega}= \cup_{K \in \calT} \overline{K}$}.
    For each control volume $K\in\calT$, $\meas{K}>0$ denotes its volume.
 \item
    $\calE$ is the set of faces such that each face $\sigma$ is a $(d-1)$-dimensional hyperplane
    with measure
    \mbox{$\meas{\sigma}>0$}.
    For each $K \in \calT$, $\calE_K$ is the subset of
    $\calE$ such that \mbox{$\partial K  = \cup_{\sigma \in\calE_K}{\sigma}$}. Furthermore, $\mathbf{x}_\sigma$ denotes
    the face evaluation points and $\n_{K,\sigma}$ the unit vector that is normal to
    $\sigma$ and outward to $K$.

  \item
    $\calP \eqbydef \lbrace \mathbf{x}_K\rbrace_{K \in \calT}$ is the set of \emph{cell centers} such that
    $\mathbf{x}_K\in K$ and $K$ is star-shaped with respect to $\mathbf{x}_K$.
  \item
    $\calV$ is the set of vertices of the grid, corresponding to the corners of the primary grid elements $\elem \in \calM$.
  \end{enumerate}
\end{definition}

For ease of exposition, we assume $d = 2$, however, this is not a limitation and the model can be readily extended to three dimensions.
In fact, our implementation is capable of handling three-dimensional systems, which is demonstrated in the section containing the
numerical examples. To distinguish geometrical quantities between the individual subdomains, the superscripts $(\cdot)^\ff$ and $(\cdot)^\porm$ are used, which are omitted if no ambiguity arises.

\subsection{Staggered-grid scheme}
\label{sub:staggered}

A staggered-grid finite volume scheme, also known as MAC scheme \cite{harlow1965a, schneider2020}, is used in the free-flow subdomain.
Scalar quantities, such as pressure and density, are defined at the cell centers $\calP$ while degrees of freedom corresponding to the components of the velocity vector are located on the primary control volumes' faces $\calE$. As opposed to collocated schemes, where all variables are defined at the same location, this resulting scheme is stable, guaranteeing oscillation-free solutions without any additional stabilization techniques \cite{versteeg2007a}.

The current implementation is restricted to grid cells that are aligned with the
coordinate axes (hyperrectangles) such that it is assumed that $\calD$ is a structured conforming quadrilateral mesh. This is a restriction of this scheme and general interface topologies between the free-flow region and the porous medium cannot be resolved by using unstructured grids.
Although extensions to this scheme, allowing general unstructured staggered schemes, are being investigated in ongoing research, these extensions are not within the scope of this paper and will not be discussed. The coupling approach presented in this paper can be seen as the basis for more general coupling concepts.

Equations to fully outline the staggered-grid discretization scheme are listed in \cite{schneider2020} and unchanged in this work. For incorporation of boundary conditions, we refer the reader to \cite{versteeg2007a}.

\subsection{Vertex-centered finite volume scheme}
\label{sub:box}

For typical cell-centered finite volume schemes, see \cite{Droniou:2014,schneider.ea:2018b} for an overview, elements and control volumes usually coincide, $\calT \equiv \calM$. On the contrary, the \BOX scheme \cite{Hackbusch89,huber2000a} introduces control volumes that do not coincide with the primary grid elements, but are defined around the grid vertices.
In order to outline their construction, let us introduce with $\calS_\vertex \subset \calS$ the set of primary grid faces adjacent to the vertex $\vertex \in \calV$, \ie for each $\gridFace \in \calS_\vertex$ it is $\vertex \in \partial \gridFace$ .
To continue, let us introduce the set of adjacent primary grid elements $\calM_\vertex$, \ie we have $\vertex \in \partial \elem$ for each $\elem \in \calM_\vertex$.
The control volume associated with the vertex $\vertex$ of the grid is then constructed by connecting the centers $\mathbf{x}_\elem$ of all adjacent elements $\elem \in \calM_\vertex$ and the centers $\mathbf{x}_\gridFace$ of all adjacent faces $\gridFace \in \calS_\vertex$.
Please note that in three-dimensional settings, the construction furthermore considers the centers of the adjacent grid edges.
For an illustration of the control volumes constructed at the interface to the free-flow domain, see \cref{fig:interactionRegion_interface}.

Each control volume is partitioned into sub-control volumes $\scv$, and with $\calT_K$ we denote the set of sub-control volumes embedded in the control volume $K \in \calT$ such that $\overline{K} = \cup_{\scv \in \calT_K} \overline{\scv}$.
Each sub-control volume is defined as the overlap of the control volume $K$ with an adjacent element $\elem \in \calM_\vertex$, that is, $\scv = K \cap \elem$.

In the following, we want to derive the discrete formulation of \labelcref{eq:darcy}.
For the sake of readability, let us consider a general balance equation, representative for \cref{eq:darcyMassBalance,eq:darcyEnergyBalance}, of the form:
\begin{equation}
    \frac{\partial \someQuantity}{\partial t} + \nabla \cdot \someFlux = q.
    \label{eq:someConservation}
\end{equation}
Here, $\someQuantity$ is the conserved quantity and $\someFlux$ is the flux term.
Integration of \cref{eq:someConservation} over a control volume $K \in \calT$, making use of the divergence theorem, yields
\begin{equation}
    \int_{K} \frac{ \partial \someQuantity}{\partial t} \, \mathrm{d}x
    + \int_{\partial K} \someFlux \cdot \n \, \mathrm{d}\Gamma
    = \int_{K} q \, \mathrm{d}x.
    \label{eq:someConservationDiscrete1}
\end{equation}
Let us now split the integral over the control volume into integrals over the sub-control volumes, and let us introduce the discrete flux $\discFlux_{K, \sigma} \approx \int_\sigma \someFlux \cdot \n \, \mathrm{d}\Gamma$ over a face $\sigma \in \calE_K$:
\begin{equation}
    \sum_{\scv \in \calT_K} \int_{\scv} \frac{ \partial \someQuantity}{\partial t} \, \mathrm{d}x
    + \sum_{\sigma \in \calE_K} \discFlux_{K, \sigma}
    = \sum_{\scv \in \calT_K }\int_{\scv} q \, \mathrm{d}x.
    \label{eq:someConservationDiscrete2}
\end{equation}
Furthermore, let us denote with $\someQuantity_\scv = |\scv|^{-1}\int_\scv \someQuantity \, \mathrm{d}x$ and $q_\scv = |\scv|^{-1}\int_\scv q \, \mathrm{d}x$ the conserved quantity and the source term evaluated for the sub-control volume $\scv$.
Assuming that the time derivative is sufficiently smooth, this yields:
\begin{equation}
    \sum_{\scv \in \calT_K} |\scv|\frac{ \partial \someQuantity_\scv}{\partial t}
    + \sum_{\sigma \in \calE_K} \discFlux_{K, \sigma}
    = \sum_{\scv \in \calT_K } |\scv| q_\scv.
    \label{eq:someConservationDiscrete3}
\end{equation}
Finally, we approximate the time derivative by a finite difference:
\begin{equation}
    \sum_{\scv \in \calT_K} |\scv| \frac{ \left( \someQuantity_\scv \right)^{t + \Delta t}
                                   - \left( \someQuantity_\scv \right)^{t}
                                 }{\Delta t}
    + \sum_{\sigma \in \calE_K} \discFlux_{K, \sigma}
    = \sum_{\scv \in \calT_K } |\scv| q_\scv,
    \label{eq:someConservationDiscrete4}
\end{equation}
where $t$ is the current time level and $\Delta t$ denotes the time step size.
In this work, we use an implicit time discretization scheme, \ie we evaluate the discrete fluxes and source terms in \eqref{eq:someConservationDiscrete4} on the time level $t + \Delta t$.
Recall that each control volume $K$ can be uniquely associated to a primary grid vertex.
As a consequence, the same also holds for the sub-control volumes $\scv \in \calT_K$.
A discrete value of $\someQuantity$ is associated to each grid vertex, which we denote with $\hat{\someQuantity}_\vertex$, $\vertex \in \calV$, and we understand $\hat{\someQuantity}_\vertex$, $\hat{\someQuantity}_\scv$ and $\hat{\someQuantity}_K$ to refer to the same discrete value.
In this work, we use the mass lumping technique proposed by \cite{huber2000a} to improve the stability of the scheme.
Thus, we use the nodal values for storage term evaluations:
\begin{equation}
    \sum_{\scv \in \calT_K} |\scv| \frac{ \left( \hat{\someQuantity}_\scv \right)^{t + \Delta t}
                                   - \left( \hat{\someQuantity}_\scv \right)^{t}
                                 }{\Delta t}
    + \sum_{\sigma \in \calE_K} \discFlux_{K, \sigma}
    = \sum_{\scv \in \calT_K } |\scv| q_\scv.
    \label{eq:someConservationDiscreteFinal}
\end{equation}
The above equation is the discrete form of \eqref{eq:someConservation} used in the \BOX scheme.
In the remainder of this section, the expressions for the discrete fluxes $\discFlux_{K, \sigma}$ are presented.
In the \BOX scheme, fluxes are computed using continuous piecewise linear basis functions on the primary discretization.
A basis function $\basisFunc_\vertex$ is associated to each vertex $\vertex \in \calV$ of the discretization, which means the basis functions are also uniquely associated to the control volumes $K \in \calT$.
Moreover, the basis function associated with the vertex $\vertex$ is equal to one at the vertex itself and zero on all other vertices of the grid:
\begin{equation}
    \basisFunc_\vertex \left(\mathbf{x}_w\right)
            = \delta_{vw}, \,\, w \in \calV,
\end{equation}
where $\mathbf{x}_w$ is the position of the vertex $w$ and $\delta_{vw}$ denotes the Kronecker-Delta.
As seen in the above condition, the basis functions have a small support region, \ie they are only non-zero within the elements adjacent to $\vertex$:
\begin{equation}
\text{supp}(\basisFunc_\vertex) = \bigcup_{\elem \in \calM_\vertex}\overline{E}.
\end{equation}
The basis functions further describe a partition of unity, and we can express the discrete approximation of $\someQuantity$ at a position $\mathbf{x} \in \Omega$ with:
\begin{equation}
    \tilde{\someQuantity} \left(\mathbf{x}\right) = \sum_{\vertex \in \calV} \hat{\someQuantity}_\vertex \basisFunc_\vertex \left(\mathbf{x}\right).
    \label{eq:discreteQuantity}
\end{equation}
Extending from \eqref{eq:discreteQuantity}, the expression for the discrete gradient of $\someQuantity$ is:
\begin{equation}
    \nabla \tilde{\someQuantity} \left(\mathbf{x}\right) = \sum_{\vertex \in \calV} \hat{\someQuantity}_\vertex \nabla \basisFunc_\vertex \left(\mathbf{x}\right).
    \label{eq:discreteQuantityGradient}
\end{equation}
Let us further introduce $\calV_\elem$, or the set of vertices contained in the primary grid element $\elem \in \calM$, which enables us to express the discrete approximation of $\someQuantity$ as well as the gradient at a position $\mathbf{x} \in \elem$ by:
\begin{equation}
    \tilde{\someQuantity} \left(\mathbf{x}\right) = \sum_{\vertex \in \calV_\elem} \hat{\someQuantity}_\vertex \basisFunc_\vertex \left(\mathbf{x}\right), \quad \quad
    \nabla \tilde{\someQuantity} \left(\mathbf{x}\right) = \sum_{\vertex \in \calV_\elem} \hat{\someQuantity}_\vertex \nabla \basisFunc_\vertex \left(\mathbf{x}\right).
\end{equation}

Making use of the above definitions, we define the discrete fluxes of the component $\kappa$ across a face $\sigma \in \calE_K$, which describes a part of the boundary of the control volume $K \in \calT$, and which is contained in the primary grid element $\elem \in \calM$, as follows:
\begin{subequations}
    \begin{align}
        \tilde{\mathbf{v}}_{\alpha, K, \sigma}
            &:= - \mathbf{K}_{\elem}
                \sum_{\vertex \in \calV_{\elem}} \hat{p}_{\alpha, \vertex} \nabla \basisFunc_\vertex (\mathbf{x}_\sigma) , \\
        \discFlux_{K, \sigma}^{\kappa, \adv}
            &:= |\sigma| \sum_\alpha \left( {\rho_\alpha X_\alpha^\kappa \lambda_\alpha} \right)_\upFactor^{\up, \sigma}
              \tilde{\mathbf{v}}_{\alpha, K, \sigma},
              \label{eq:discreteMassFluxesAdv} \\
        \discFlux_{K, \sigma}^{\kappa, \diff}
            &:= |\sigma| \sum_\alpha \rho_{\alpha, \sigma} D_{\alpha, \sigma}
              \sum_{\vertex \in \calV_{\elem}} \hat{X}^\kappa_{\alpha, \vertex} \nabla \basisFunc_\vertex (\mathbf{x}_\sigma).
    \label{eq:discreteMassFluxesDiff}
    \end{align}
\label{eq:discreteMassFluxes}
\end{subequations}
Here, we introduced with $\mathbf{K}_{\elem}$ the permeability defined for the primary grid element ${\elem}$, and with $\rho_{\alpha, \sigma}$ and $D_{\alpha, \sigma}$ the phase density and diffusion coefficient associated with the face $\sigma$, respectively.
While the density at the face is defined by the interpolation of the nodal values $\rho_{\alpha, \sigma} = \sum_{\vertex \in \calV_{\elem}} \hat{\rho}_{\alpha, \vertex} \basisFunc_\vertex (\mathbf{x}_\sigma) $, the diffusion coefficient on the face is taken as the harmonic mean of the diffusion coefficients in the two adjacent sub-control volumes.

Similarly, we define the discrete heat fluxes as:
\begin{subequations}
    \begin{align}
        \discFlux_{K, \sigma}^{h, \adv}
            &:= |\sigma|\sum_\alpha \left( \rho_\alpha h_\alpha \right)_\upFactor^{\up, \sigma}
              \tilde{\mathbf{v}}_{\alpha, K, \sigma},
              \label{eq:discreteHeatFluxesConv}
              \\
        \discFlux_{K, \sigma}^{\cond}
            &:= -|\sigma| \lambda_{\mathrm{eff}, \sigma}
               \sum_{\vertex \in \calV_{\elem}} \hat{T}_{\vertex} \nabla \basisFunc_\vertex (\mathbf{x}_\sigma).
    \label{eq:discreteHeatFluxesCond}
    \end{align}
\label{eq:discreteHeatFluxes}
\end{subequations}
Here, the effective thermal conductivity on the face is again defined as the harmonic mean of the conductivities in the two sub-control volumes adjacent to it.
In the discrete fluxes stated in \eqref{eq:discreteMassFluxesAdv} and \eqref{eq:discreteHeatFluxesConv} we have used the notation $\left( u \right)_\upFactor^{\up, \sigma}:= \upFactor   u^\up + (1-\upFactor) u^\dn$ to denote quantities that are evaluated in the upstream or downstream sub-control volume of face $\sigma$, i.e. $\upFactor = 1$ (used as default value) corresponds to a full upwinding and $\upFactor = 0.5$ to an averaging.
With the discrete conservation \cref{eq:someConservationDiscreteFinal} and the discrete fluxes \eqref{eq:discreteMassFluxes} and \eqref{eq:discreteHeatFluxes}, we can state the discrete form of the mass balance \cref{eq:darcyMassBalance} for each control volume $K \in \calT$:
\begin{equation}
    \phi \sum_{\scv \in \calT_K} |\scv| \frac{ \left(
                                          \hat{\rho}_{\alpha, \scv}
                                          \hat{X}_{\alpha, \scv}^\kappa
                                          \hat{S}_{\alpha, \scv}
                                   \right)^{t + \Delta t}
                                   - \left(
                                          \hat{\rho}_{\alpha, \scv}
                                          \hat{X}_{\alpha, \scv}^\kappa
                                          \hat{S}_{\alpha, \scv}
                                     \right)^{t}
                                 }{\Delta t}
    + \sum_{\sigma \in \calE_K} \discFlux_{K, \sigma}^{\kappa, \adv} + \discFlux_{K, \sigma}^{\kappa, \diff}
    = \sum_{\scv \in \calT_K }|\scv| q_\scv,
    \label{eq:discreteDarcyMass}
\end{equation}
and correspondingly for the energy balance \cref{eq:darcyEnergyBalance}:
\begin{equation}
    \sum_{\scv \in \calT_K} |\scv| \frac{ \left( \hat{e} \right)^{t + \Delta t}
                                   - \left( \hat{e} \right)^{t}
                                 }{\Delta t}
    + \sum_{\sigma \in \calE_K} \discFlux_{K, \sigma}^{h, \adv} + \discFlux_{K, \sigma}^{\cond}
    = \sum_{\scv \in \calT_K } |\scv| q_\scv,
    \label{eq:discreteDarcyEnergy}
\end{equation}
where $\hat{e} = \left(1-\phi\right) \hat{\rho}_s \hat{c}_s \hat{T}_s + \phi \sum_\alpha \hat{\rho}_{\alpha, \scv} \hat{u}_{\alpha, \scv} \hat{S}_{\alpha, \scv}$.


\subsection{Coupling}
\label{sub:discrete_coupling}
\begin{figure}[ht!]
  \centering
  \includegraphics[width=1.0\linewidth,keepaspectratio]{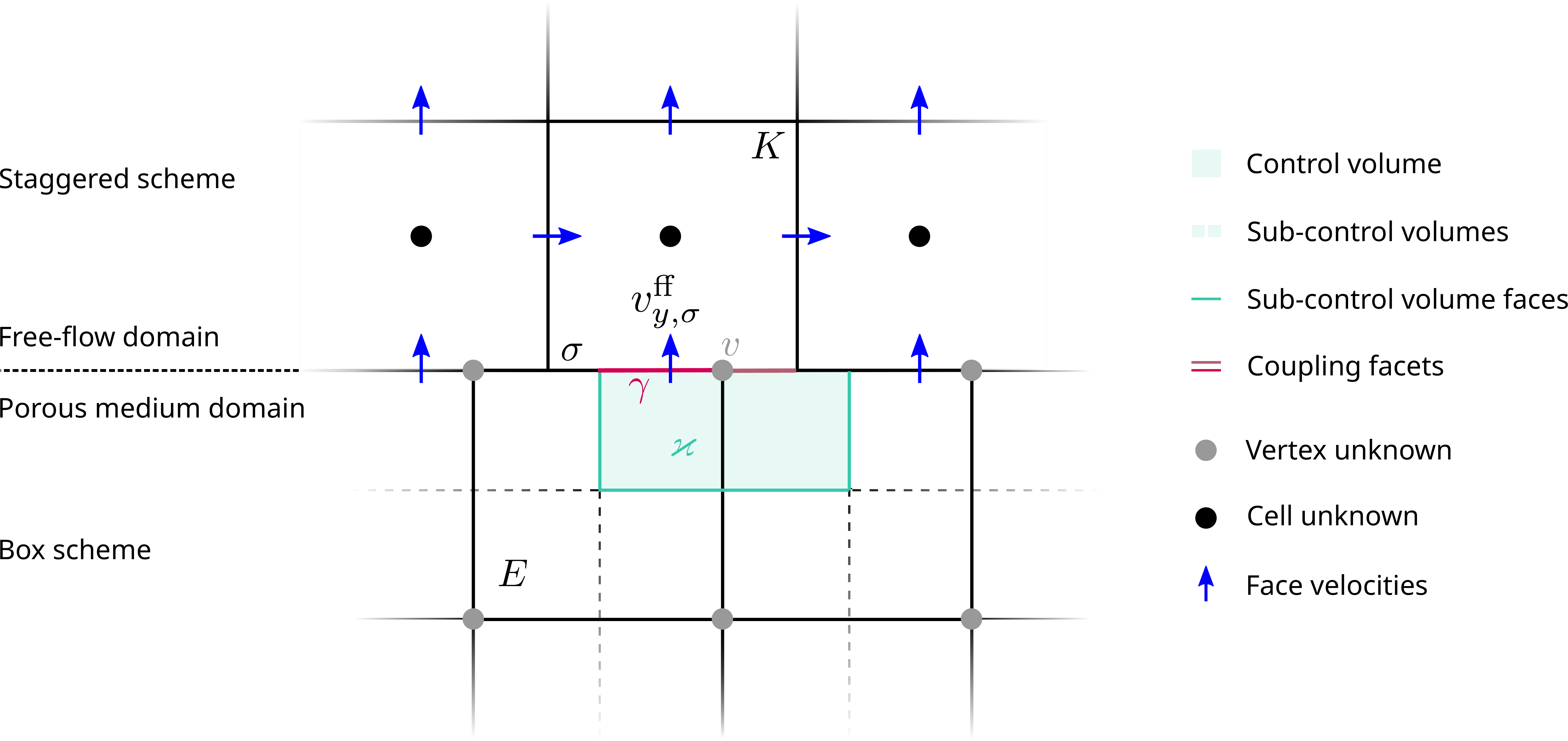}
  \caption{Control volume of the \BOX method at the interface to the free-flow subdomain. The unknowns $u_v^\porm$ of the \BOX method are located at grid vertices, whereas pressures, mass fractions, and temperatures of the staggered scheme are located at cell centers, and velocities at faces. The basis functions of the \BOX method can be used to interpolate solution values at the interface.}
  \label{fig:interactionRegion_interface}
\end{figure}

This section is devoted to the incorporation of the coupling conditions \labelcref{eq:interfaceConditions} in the discrete setting.
The discretizations of the sub-domains are allowed to be non-matching at the interface, which results in arbitrary overlap between grid faces of the free-flow and the porous-medium domain, respectively.
Let us partition the interface $\Gamma^\ipmff$ into a set of facets $\iSectionSet$  (line segments for the two-dimensional case), for which we note that $\overline{\Gamma}^\ipmff = \bigcup_{\iSection \in \iSectionSet} \overline{\iSection}$.
Each facet $\iSection \in \iSectionSet$ is the intersection between a face of the discretization of the porous-medium domain and a face of the discretization of the free-flow domain, depicted as coupling segments in \cref{fig:interactionRegion_interface}.
Each face $\sigma$ of the discretizations of the subdomains overlaps with one or more such facets, which we collect in the set $\iSectionSet_\sigma$ and where we note that $\overline{\sigma} = \cup_{\iSection \in \iSectionSet_\sigma} \overline{\iSection}$.

The coupling between the two subdomains is realized as follows: \cref{eq:ifMassFluxCont,eq:ifTempFluxCont} represent Neumann (flux) boundary conditions, where only the free-flow terms are calculated and set accordingly as boundary fluxes. This guarantees the conservation of mass and energy.
To calculate these fluxes, we use the advantage of the vertex-centered finite volume scheme, where degrees of freedom are naturally located at the interface.
In the following, we present the construction of the free-flow fluxes at the interface by using a two-point flux approximation (Tpfa) \cite{schneider.ea:2018b} to discretize diffusive and conductive fluxes given by Fick's and Fourier's law \eqref{eq:FicksLaw}-\eqref{eq:FouriersLaw}. This choice is based on the fact that, in the free-flow domain, we are restricted to grid cells that are aligned with the coordinate axes, as mentioned in \cref{sub:staggered}. In this case, the Tpfa yields an accurate flux approximation.

Let $\sigma \in \calE^\ff$ be a face of some control volume $K \in \calT^\ff$ of the free-flow grid (see \cref{fig:interactionRegion_interface}) such that $\sigma \subset \Gamma^\ipmff$. Integrating the free-flow related terms of \cref{eq:ifMassFluxCont,eq:ifTempFluxCont} over $\sigma$ yields the following discrete approximations:
\begin{align}
 \int_\sigma \left[\left(\varrho_g X_g^\kappa \vel_g + \mathbf{J}^{\kappa, \diff} \right)\cdot \n \right]^\ff \mathrm{d}\Gamma
 &\approx \sum_{\iSection \in \iSectionSet_\sigma} \left[  \discFlux_{K, \iSection}^{\kappa, \adv} +  \discFlux_{K, \iSection}^{\kappa, \diff} \right], \\
 \int_\sigma  \left[ \left(\varrho_g h_g \vel_g + \sum_{\kappa} h^{\kappa}_g  \mathbf{J}^{\kappa, \diff} + \mathbf{J}^\cond \right)\cdot \n \right]^\ff  \mathrm{d}\Gamma
 &\approx \sum_{\iSection \in \iSectionSet_\sigma}  \left[ \discFlux_{K, \iSection}^{h,\adv} + \sum_\kappa  \discFlux_{K, \iSection}^{\kappa, \diff} + \discFlux_{K, \iSection}^{\cond} \right], \end{align}
with flux approximations
\begin{align}
\label{eq:ifadvfluxes}
  \discFlux_{K, \iSection}^{\kappa, \adv} := |\iSection|(\varrho_g X_g^\kappa)_\upFactor^{\up,\iSection} \vel^\ff_{g,\sigma}, \quad &
  \discFlux_{K, \iSection}^{h,\adv} := |\iSection|(\varrho_g h_g)_\upFactor^{\up,\iSection}  \vel^\ff_{g,\sigma} \cdot \mathbf{n},
  \\
  \label{eq:ifdifffluxes}
 \discFlux_{K, \iSection}^{\kappa, \diff} := |\iSection|\varrho_{g,\iSection} D_{g ,\iSection}^{\kappa} \frac{(X_{g,K}^{\kappa ,\ff} - \Pi_\sigma \tilde{X}_{g}^{\kappa,\porm})}{d_{K,\sigma}}, \quad &  \discFlux_{K, \iSection}^{\cond} := |\iSection| \lambda_{\text{eff},\iSection} \frac{(T_{K}^{\ff} - \Pi_\sigma \tilde{T}^\porm)}{d_{K,\sigma}},
\end{align}
which are related to the control volume $K\in \calT^\ff$ and the corresponding sub-control volume $\scv$ of the porous-medium domain such that $\iSection = \overline{K}\cap \overline{\scv}$. The interface density is evaluated with the following average, $\varrho_{g,\iSection} := 0.5(\varrho_{g,K}^\ff + \varrho_{g,\scv}^\porm)$ but the diffusion coefficient is evaluated as $D_{g ,\iSection}^{\kappa} := D_{g,K}^{\kappa,\ff}$, using only information from the free-flow control volume. The upwind quantities are analogous to the previous definition and the default is $\upFactor = 1$ such that only upstream information, i.e. either from the free-flow or the the porous-medium domain, is used.

For the construction of diffusive/conductive fluxes \labelcref{eq:ifdifffluxes} we use the continuity of the primary variables $X_{g}^{\kappa}$ and $T$, which is provided via conditions \labelcref{eq:ifMassFracCont,eq:ifTempCont}. This allows an evaluation of these quantities using information from the discrete porous-medium solution.
This is done with two types of projection operators $\Pi_\sigma \colon H^1(\Omega^\porm) \mapsto \R$:
\begin{align}
\label{eq:l2projection}
\ltwoproj_\sigma \tilde{\someQuantity} &:= \frac{1}{|\sigma|} \int_\sigma \tilde{\someQuantity} \, \mathrm{d}\Gamma = \frac{1}{|\sigma|} \sum_{{\iSection \in \iSectionSet_\sigma}} \int_\iSection \tilde{\someQuantity} \, \mathrm{d}\Gamma, \\
\label{eq:avgprojection}
\avgproj_\sigma \tilde{\someQuantity} &:=\frac{1}{|\sigma|} \sum_{{\iSection \in \iSectionSet_\sigma}, v \in \overline{\iSection}}  |\iSection| \hat{\someQuantity}_v.
\end{align}
$\avgproj_\sigma$ can be seen as a simplification of \labelcref{eq:l2projection} by using a simple quadrature rule, which has the advantage of yielding smaller coupling stencils. By using the same fluxes also in the porous domain, i.e. $\discFlux^\porm_{\scv, \iSection} = - \discFlux^\ff_{K, \iSection}$, the conservation of advective/convective and diffusive/conductive fluxes can directly be observed.

The remaining coupling conditions \eqref{eq:ifNormalMomentum}, for normal momentum transfer, and the slip condition \eqref{eq:ifBJS}, are needed only for the free-flow subdomain. \Cref{eq:ifNormalMomentum} yields a Neumann boundary condition for the momentum balance by evaluating $p_g^\porm$ at each free-flow face. This is again realized by using the projections $\Pi_\sigma$ introduced above. When using the Saffman simplification of the Beavers-Joseph condition, as shown in \Cref{eq:ifBJS}, the porous-medium velocity is neglected meaning only free-flow unknowns are required and the condition is independent of the chosen porous-medium discretization scheme.

\section{Numerical results}
\label{sec:numresults}

The coupling approach presented in the previous section was implemented
within the the open-source simulation environment \Dumux \cite{dumux}, which exists as an additional module of the DUNE project \citep{blatt.ea:2016}.
We employ a monolithic approach, where both sub-problems are assembled into one system of equations and use an implicit Euler method for the time discretization.
To accommodate for the nonlinearities in the systems of equations, Newton's method is employed
and the direct solver UMFPack \citep{Davis:2004a} is used to solve the resulting linear systems of equations produced in each Newton iteration. The Dune-Subgrid \cite{Graeser2009a} and Dune-UGGrid modules are used to represent the structured grids in the free-flow
and the unstructured grids used in the porous-medium subdomain, respectively.

Within this section, a series of three test cases will be used to examine the performance of the proposed coupling approach.
\Cref{subsec:convtest} provides a convergence study, while \cref{sec:fetzer} presents a comparison between the proposed method
and the coupling approaches introduced in \cite{fetzer2017a}. Finally, a showcase simulation is given in \cref{subsec:filtertest} to demonstrate the flexibility of the method at application scale.

\subsection{Convergence Study}
\label{subsec:convtest}
In this section we investigate the convergence behavior of the proposed coupling scheme on a test case that was
originally presented in~\cite{schneider2020}. It assumes incompressible single-phase flow, neglecting non-isothermal,
compositional, and gravitational effects. Consider a domain $\Omega = (0, 1) \times (0, 2)$, decomposed into
$\Omega^\porm = (0, 1)^2$ and $\Omega^\ff = (0, 1) \times (1,2)$ with interface $\Gamma^\ipmff = \{(x,y)^T \in \Omega \,|\, y = 1\}$. The material parameters are
given by
\begin{equation}
	\mu = \varrho = \BFCoeff = 1, \quad \quad
	\mathbf{K} = \begin{bmatrix}
		1 & -\frac{c}{2\omega} \sin(\omega x) \\ -\frac{c}{2\omega} \sin(\omega x) & e^{-2}(1+c\cos(\omega x))
	\end{bmatrix},
\end{equation}
where, depending on the parameter $c$, the permeability tensor depends on the spatial coordinate $x$ and varies with
the angular frequency $\omega$. Furthermore, let the solution be given by
\begin{align}
	\bm{v}^\ff &:=
	\begin{bmatrix}
		y \\
		- y \sin(\omega x)
	\end{bmatrix}, &
		\bm{v}^\porm &:=
	\begin{bmatrix}
		\frac{c}{2\omega}e^{y + 1}\sin^2(\omega x)  - \omega (e^{y + 1} + 2 - e^2) \cos(\omega x) \\
		\left(\frac{c}{2}\cos(\omega x)(e^{y + 1} + 2 - e^2) - (1 + c \cos(\omega x))e^{y - 1}\right) \sin(\omega x)
	\end{bmatrix}, \\
	p^\ff &:= -y^2 \sin^2(\omega x),  &
	p^\porm &:= (e^{y + 1} + 2 - e^2) \sin(\omega x).
\end{align}
Neglecting time derivatives, we have in the free-flow domain
\begin{align}
  \nabla \cdot \vel^\ff &= q^\ff:= -\sin(\omega x), \\
  \nabla \cdot \left(\vel \vel^{\mathrm{T}}
  - (\nabla \vel + (\nabla \vel)^{\mathrm{T}})
  + p \mathbf{I} \right)^\ff
  &= \textbf{f}
  := \begin{bmatrix}
	-2 \omega y^2 \sin(\omega x) \cos(\omega x) - 2 y \sin(\omega x) + \omega \cos(\omega x) \\
	-\omega y^2 \cos(\omega x) - \omega^2 y \sin(\omega x)
  \end{bmatrix},
\end{align}
and in the porous medium
\begin{align}
  \nabla \cdot \vel^\porm &= q^\porm := \left( 1.5c e^{y + 1}\cos(\omega x) + \omega^2 (e^{y + 1} + 2 - e^2) -(1 + c\cos(\omega x))e^{y - 1} \right) \sin(\omega x), \\
  \vel^\porm + \mathbf{K} \nabla p^\porm &= 0.
\end{align}
The validity of the coupling conditions can also be easily verified.
In order to guarantee that $\mathbf{K}$ is positive definite and a full tensor with
non-negligible off-diagonal entries (depending on $x$), we set $\omega = \pi$ and $c=0.9$.
For calculating the fluxes in the porous domain, the permeability tensors are evaluated at the cell
centers. The functions $q^\porm, q^\ff$ , and $\mathbf{f}$ are numerically integrated using a 5th-order
quadrature rule.

We use the following discrete pressure error norm,
\begin{equation}
  \error_{p_i}^m = \left(\sum_{K \in \calT_m} \meas{K} \left( p_K - \overline{p}_K \right)^2 \right)^\frac{1}{2},
  \label{eq:discreteNormPressure}
\end{equation}
where $i \in \{\porm, \ff\}$ and where the parameter $m$ indicates the $m$-th refinement level.
Furthermore, $p_K$ and $\overline{p}_K$ denote the
control volume-wise constant discrete and exact pressure values, obtained by evaluating both in the center of the
control volumes $K \in \calT_m$.

For the staggered-grid scheme, we additionally quantify the errors for the velocity unknowns. These errors are calculated as
\begin{equation}
\error_{v_i}^m = \left(\sum_{K_i^* \in \calT_{i,m}^*} \meas{K_i^*} \left( v_{i,K_i^*} - \overline{v}_{i,K_i^*} \right)^2 \right)^\frac{1}{2}, \quad i \in \{ x, y \},
  \label{eq:discreteNormVelocity}
\end{equation}
where the velocity unknowns $v_{i,K_i^*}$ are associated with the dual control volumes $K_i^*,  i \in \{ x, y \}$,
which are constructed around the faces $\sigma \in \mathcal{E}$ (see \cref{sec:discretization}). The exact velocity
values at the face centers are denoted as $\overline{v}_{i,K_i^*}$.

\begin{figure}
	\centering
	\begin{subfigure}{0.2\textwidth}
		\centering
		\includegraphics[width=\textwidth]{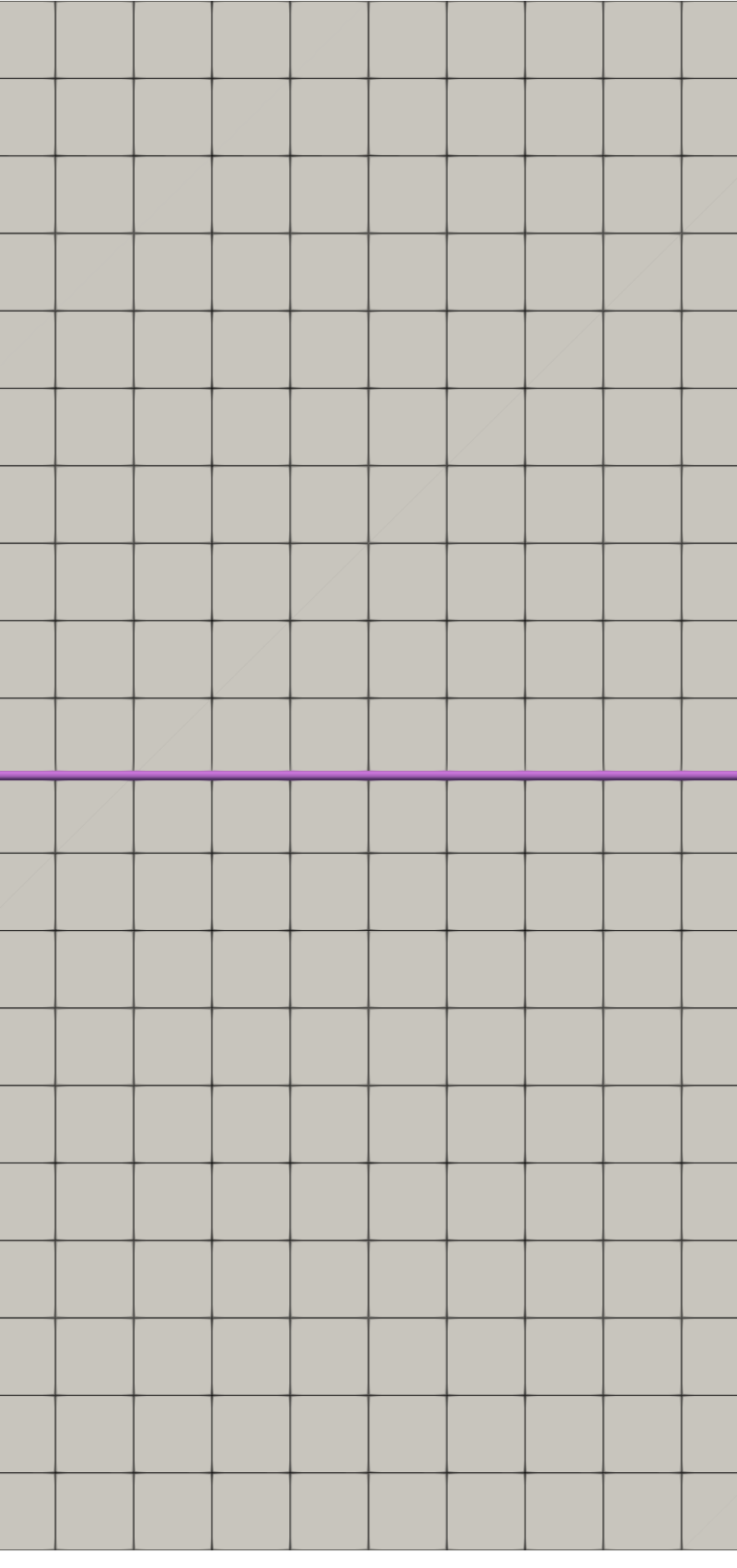}
		\caption{}
		\label{fig:convtest_grids_conforming}
	\end{subfigure}
	\begin{subfigure}{0.2\textwidth}
		\centering
		\includegraphics[width=\textwidth]{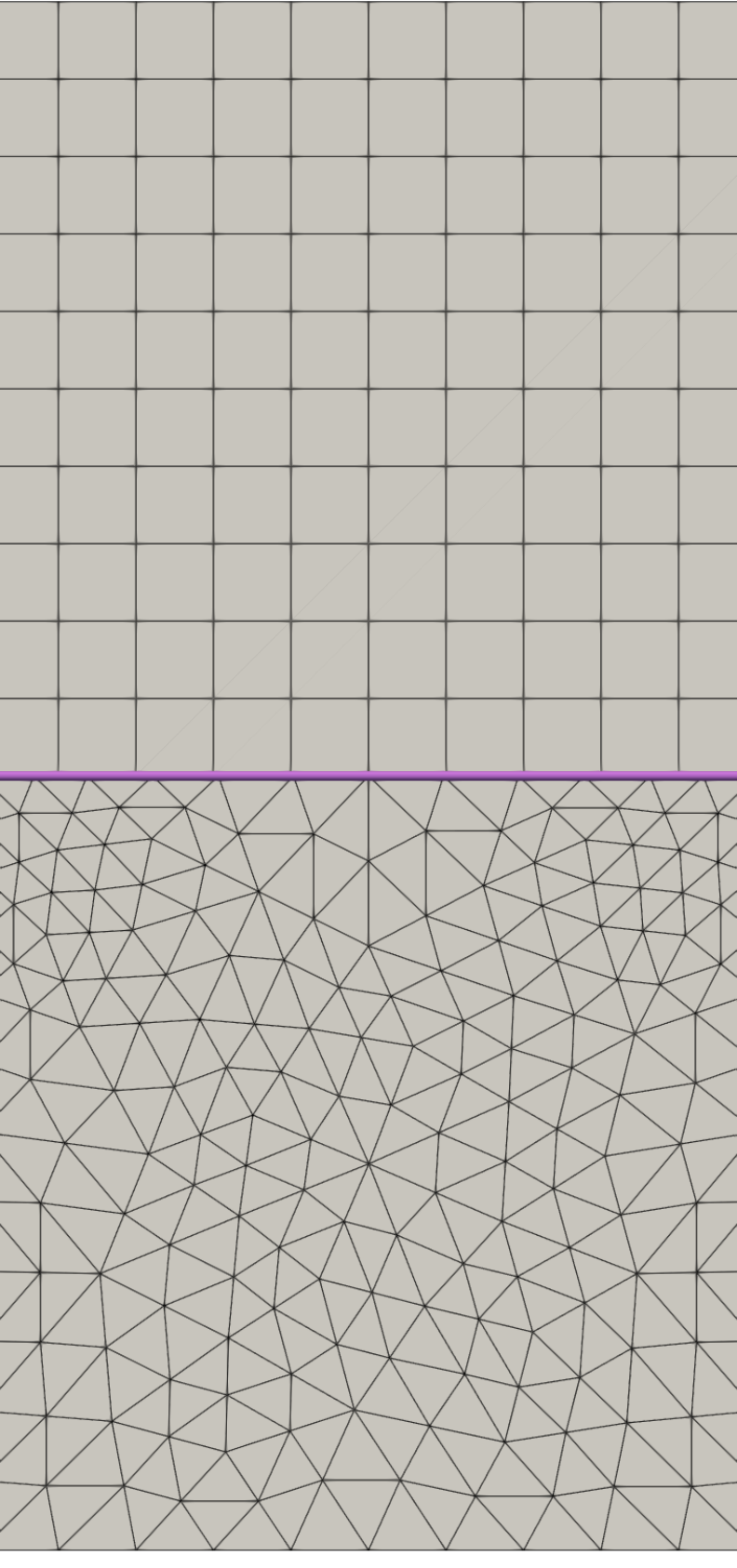}
		\caption{}
		\label{fig:convtest_grids_nonconforming}
	\end{subfigure}
	\begin{subfigure}{0.2\textwidth}
		\centering
		\includegraphics[width=\textwidth]{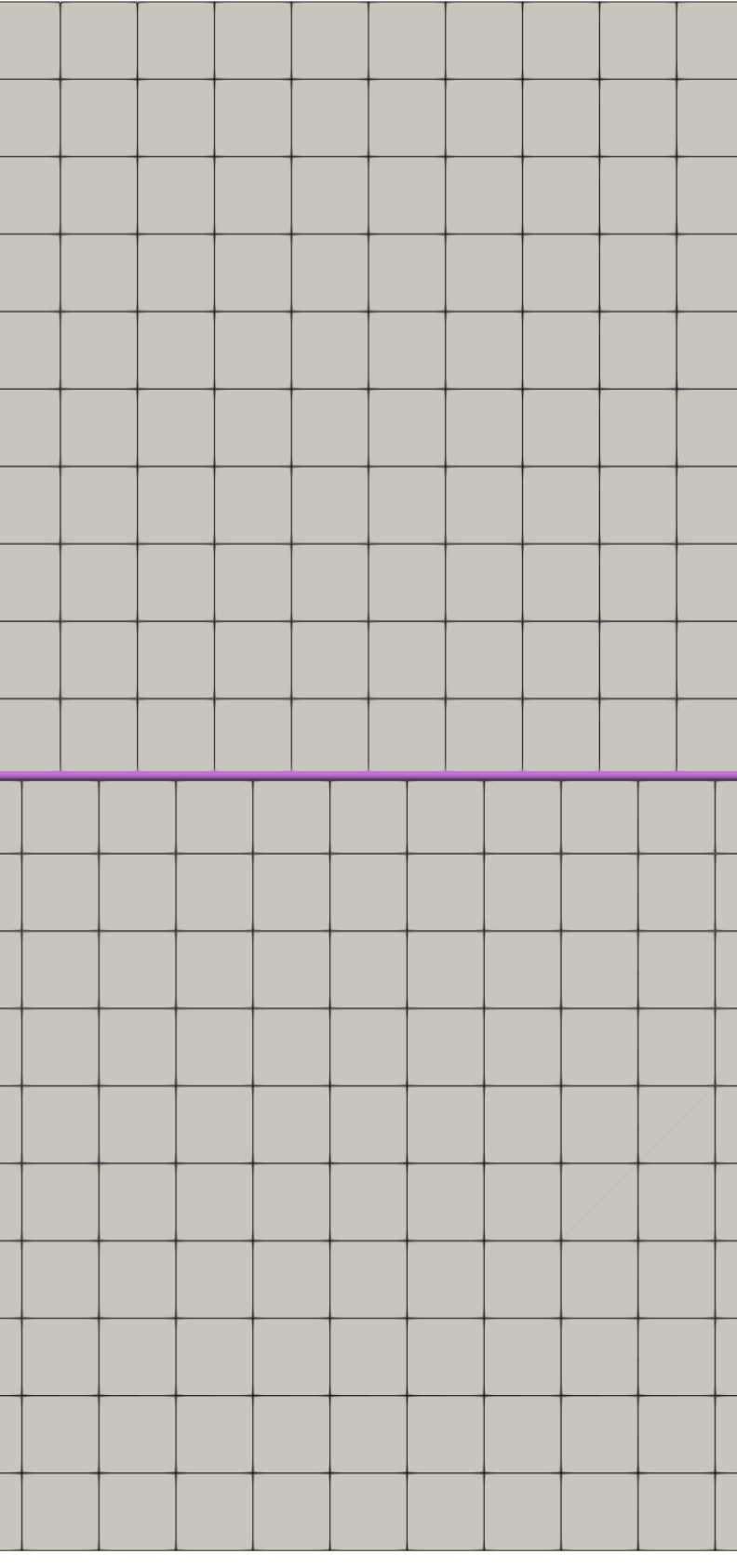}
		\caption{}
		\label{fig:convtest_grids_boxconforming}
	\end{subfigure}
    \caption{Three grid configurations used in the convergence test: (\subref{fig:convtest_grids_conforming}) conforming, (\subref{fig:convtest_grids_nonconforming}) non-conforming, and (\subref{fig:convtest_grids_boxconforming}) box-conforming. All grids are shown after the first refinement.}
	\label{fig:convtest_grids}
\end{figure}

Three different grid configurations in the porous-medium domain are considered: a grid that is conforming to the grid used
in the free-flow domain (\subref{fig:convtest_grids_conforming}), an unstructured simplex grid with non-matching interface (\subref{fig:convtest_grids_nonconforming}), and a structured grid for which the
dual grid used in the \BOX scheme is conforming to the grid used in the  free-flow domain (\subref{fig:convtest_grids_boxconforming}).
These cases are refered to here as the conforming, non-conforming, and box-conforming cases, respectively.
\Cref{fig:convtest_grids} shows
all configurations after the first grid refinement. Moreover, we consider both projection operators
\ltwoproj (see \cref{eq:l2projection}) and \avgproj (see \cref{eq:avgprojection}).

\begin{figure}
    \begin{subfigure}{0.49\textwidth}
        \centering
        \includegraphics[width=0.99\textwidth]{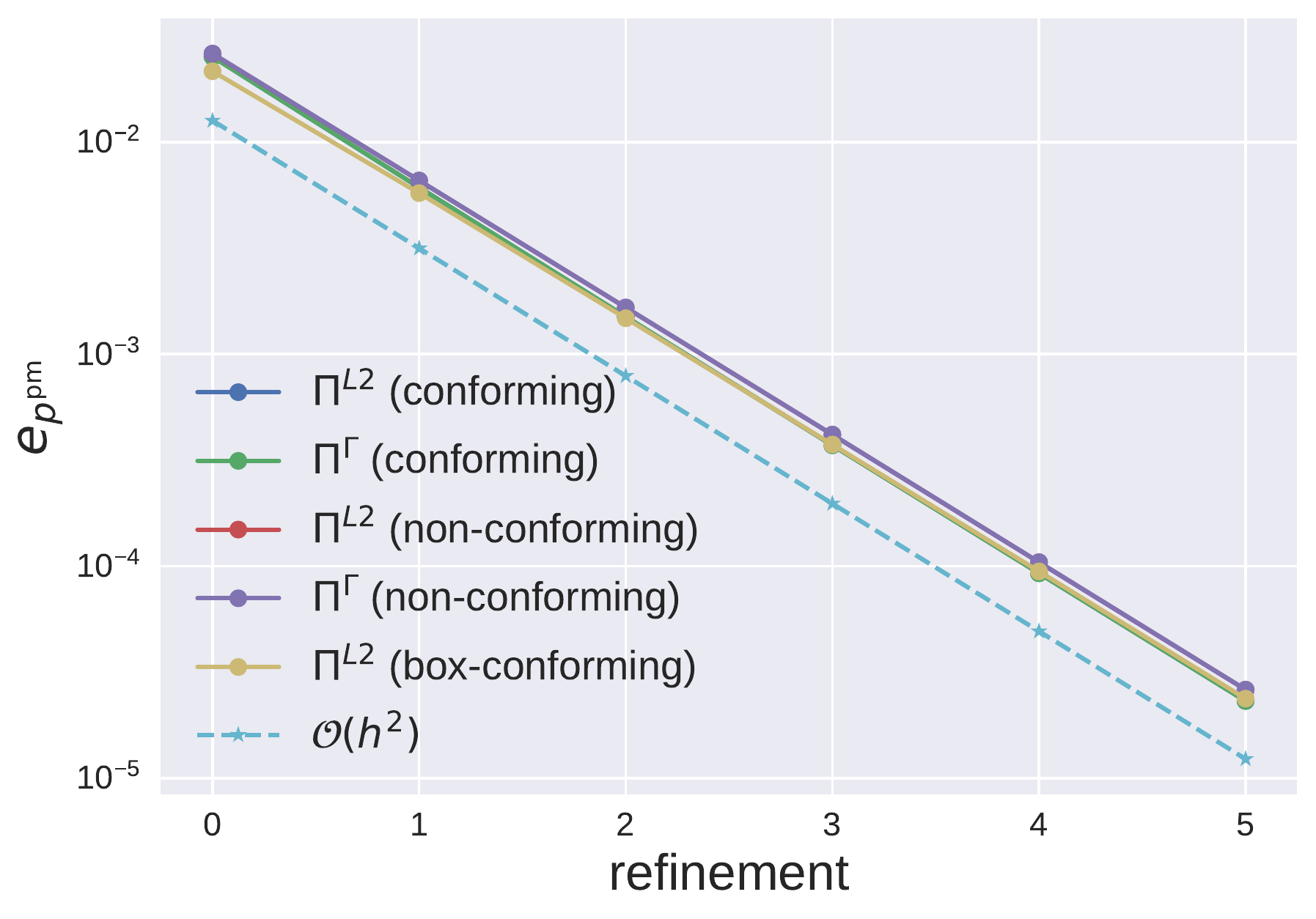}
        \caption{}
        \label{fig:convtest_darcy_p}
    \end{subfigure}
    \begin{subfigure}{0.49\textwidth}
        \centering
        \includegraphics[width=0.99\textwidth]{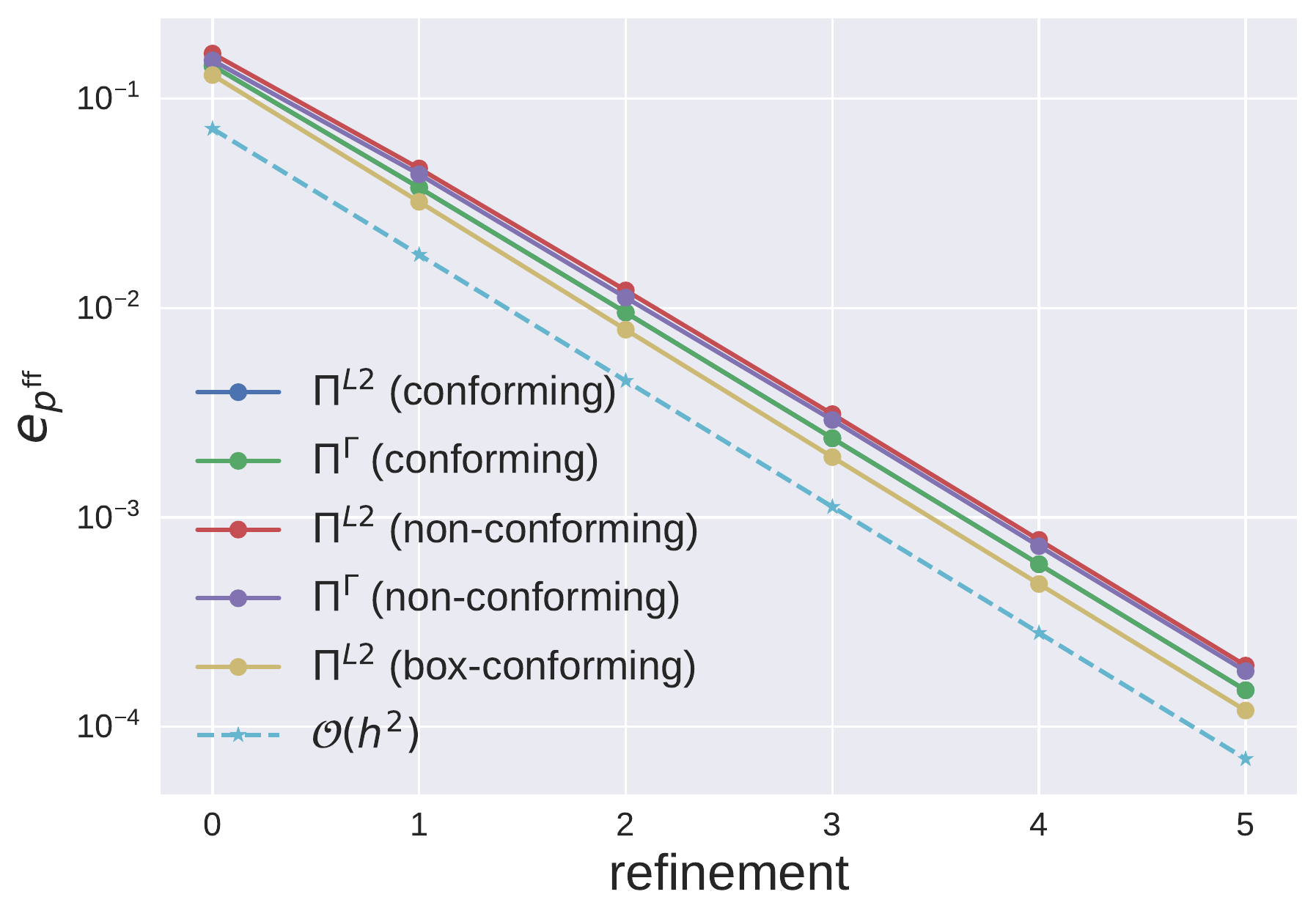}
        \caption{}
        \label{fig:convtest_stokes_p}
    \end{subfigure}
    \begin{subfigure}{0.49\textwidth}
        \centering
        \includegraphics[width=0.99\textwidth]{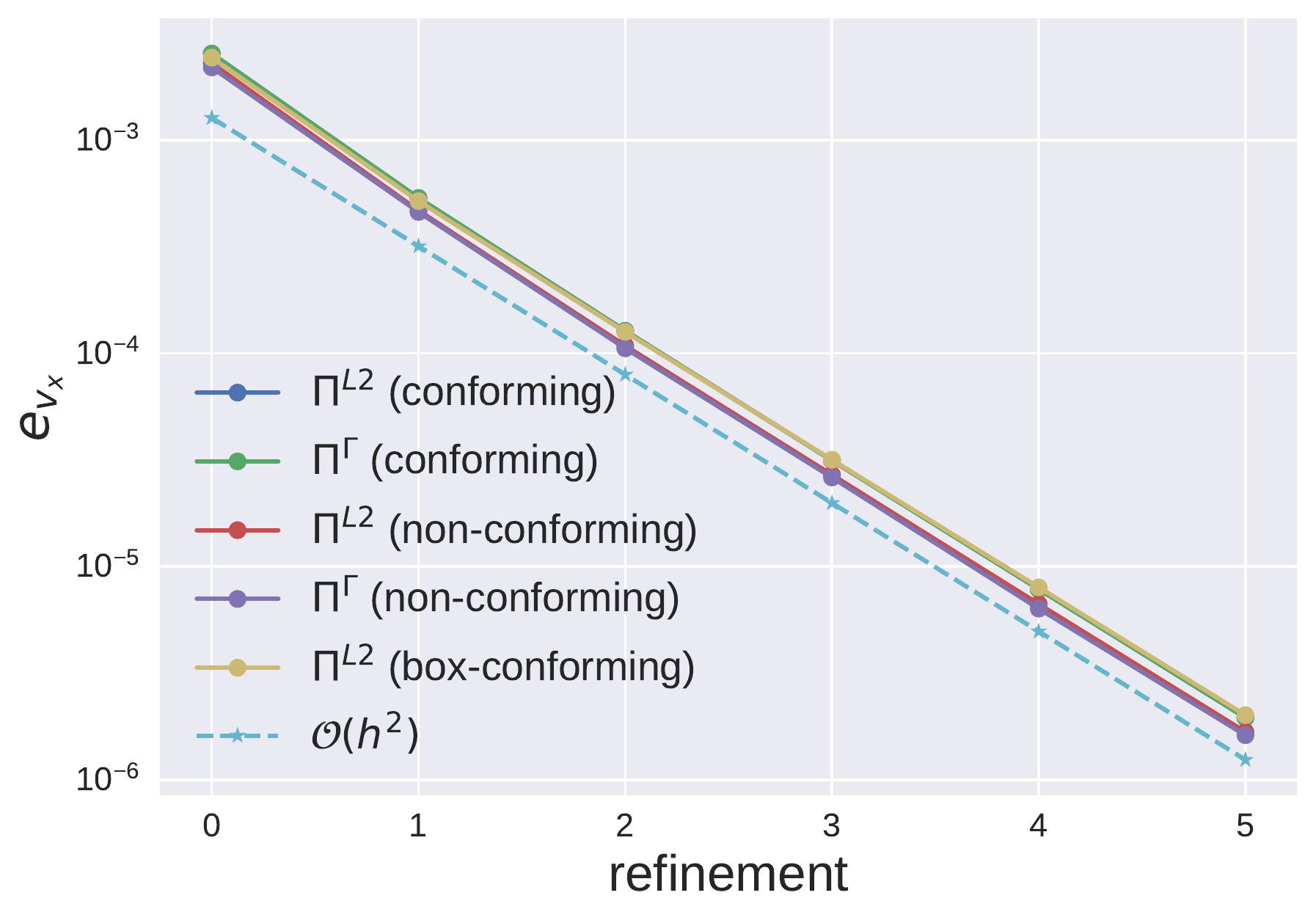}
        \caption{}
        \label{fig:convtest_stokes_vx}
    \end{subfigure}
    \begin{subfigure}{0.49\textwidth}
        \centering
        \includegraphics[width=0.99\textwidth]{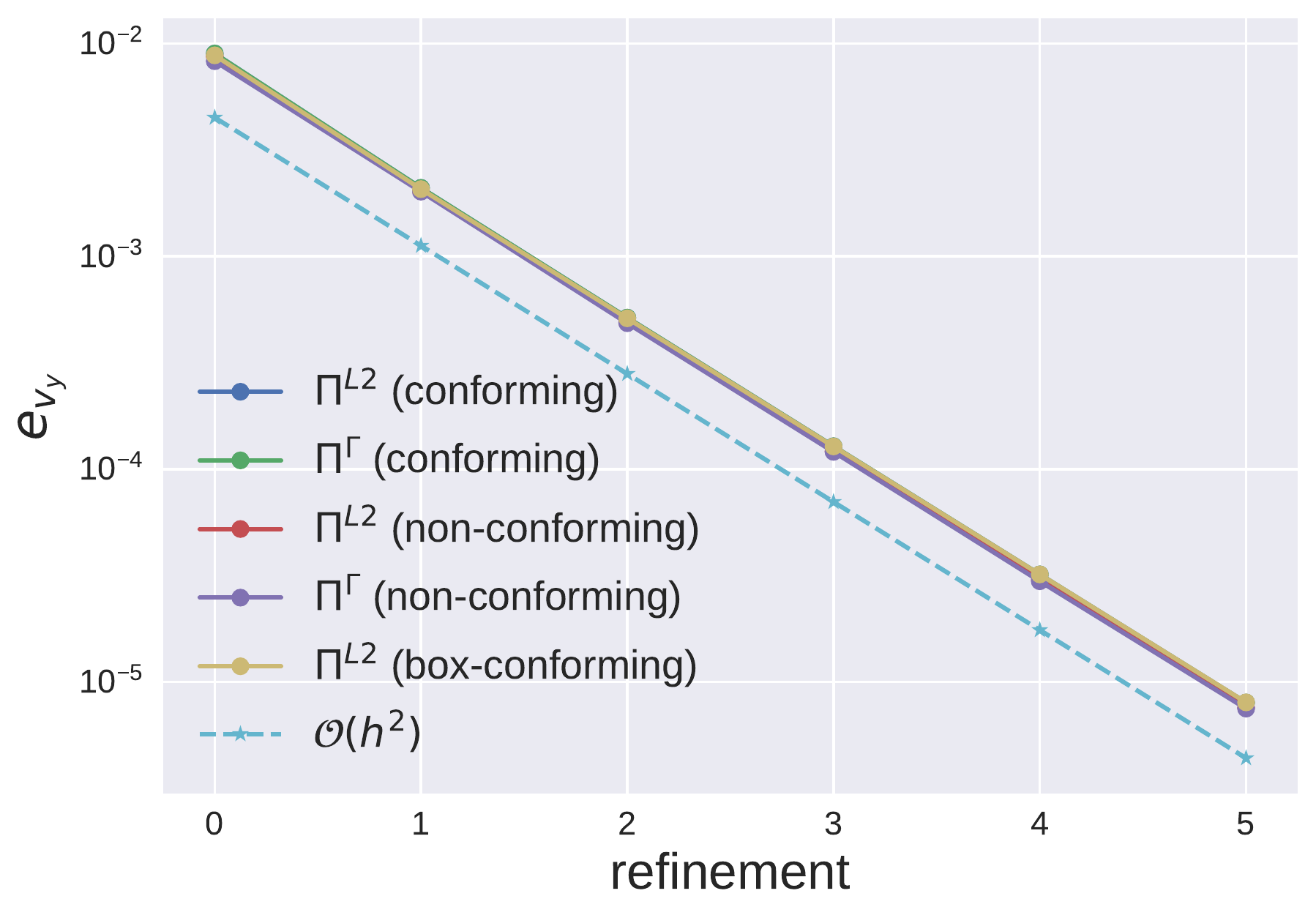}
        \caption{}
        \label{fig:convtest_stokes_vy}
    \end{subfigure}
    \caption{\textbf{Convergence test}. Error norms vs grid refinement.}
    \label{fig:convtest_plots}
\end{figure}

\Cref{fig:convtest_plots} shows the error norms plotted over grid refinement for all combinations
of grids and projection operators. Second order convergence of the considered discrete error norms
is observed in all combinations, with only minor differences in the error values. The largest differences
appear for $e_{p^\ff}$, where the smallest errors are observed with $\ltwoproj$ on the box-conforming
grid configuration, and the largest errors occur on the non-conforming grid configuration.
Note that in this case, the errors observed with the projection operator \avgproj are slightly smaller
than with \ltwoproj, despite being computationally cheaper. Note also that on the box-conforming
grid, the projections \ltwoproj and \avgproj are equivalent, which is why we only show the result for
\ltwoproj here.

\subsection{Porous obstacle} \label{sec:fetzer}
The purpose of this test is to show that the presented coupling approach is able to accurately
capture the interface exchange processes described by the continuous coupling conditions
\cref{eq:interfaceConditions}.
By comparing different coupling methods, \cite{fetzer2017a} have demonstrated that for an accurate
description of the processes at the interface, additional interface unknowns have to be introduced.
In their publication, a cell-centered Tpfa scheme on a structured grid has been used in the porous-medium domain coupled to a  staggered-grid finite volume scheme used in the free-flow domain.
When solving physically complex flow processes, described by nonlinear PDEs, the additional interface unknowns have to be calculated locally within a sub-routine that requires a nonlinear solver.
As mentioned before, the advantage of the presented approach is that no interface solver is needed because
the numerical schemes allow for an evaluation of all required quantities directly at the interface.

\begin{figure}[ht!]
  \centering
  \includegraphics[width=0.5\linewidth,keepaspectratio]{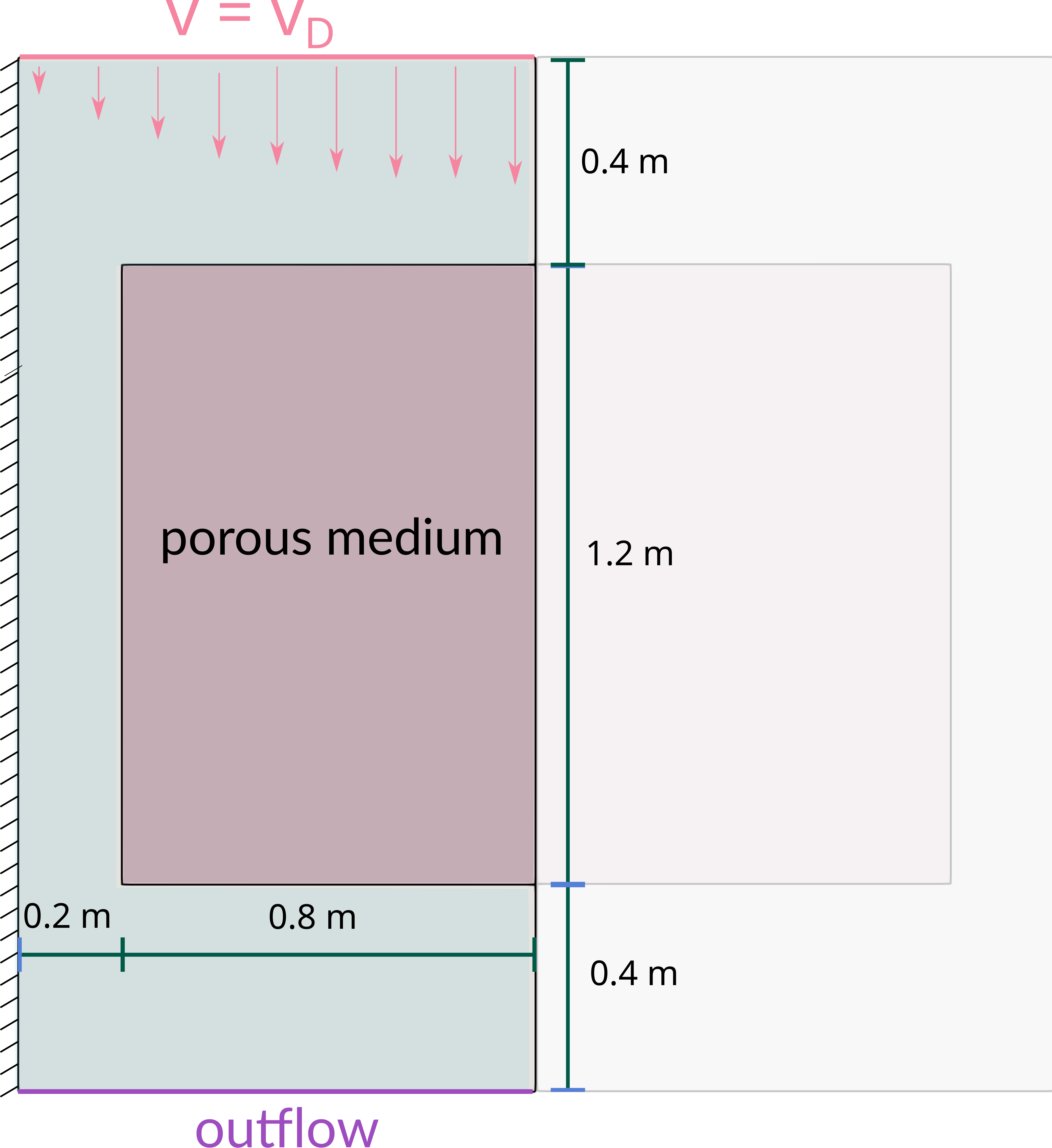}
  \caption{\textbf{Porous obstacle}. Setup and boundary conditions for the test case taken
  from \cite{fetzer2017a}.}
  \label{fig:vertical_test_setup}
\end{figure}

In this section, it is shown that the presented approach yields results that are comparable to those
presented in \cite{fetzer2017a} with the coupling method ``CM4'', which involves local nonlinear solves
for the interface quantities. In order to draw a fair comparison, we replicate a test case that was used
in the original study and which is illustrated in \cref{fig:vertical_test_setup}. It considers a
two-dimensional free-flow channel in which a flow of air is enforced by setting a parabolic
velocity profile at the top boundary, while allowing outflow across the bottom boundary. A porous medium
is placed in the center of the channel, leaving narrower channels on the left and right for the air to
potentially bypass the porous medium. Due to the symmetry of the setup, the test case only considers
the left half of the domain while applying symmetry boundary conditions on the right boundary. The
left boundary is set to no-flow. For indications on the domain dimensions, see \cref{fig:vertical_test_setup}.

\begin{figure}
    \begin{subfigure}{0.33\textwidth}
        \centering
        \includegraphics[width=0.99\textwidth]{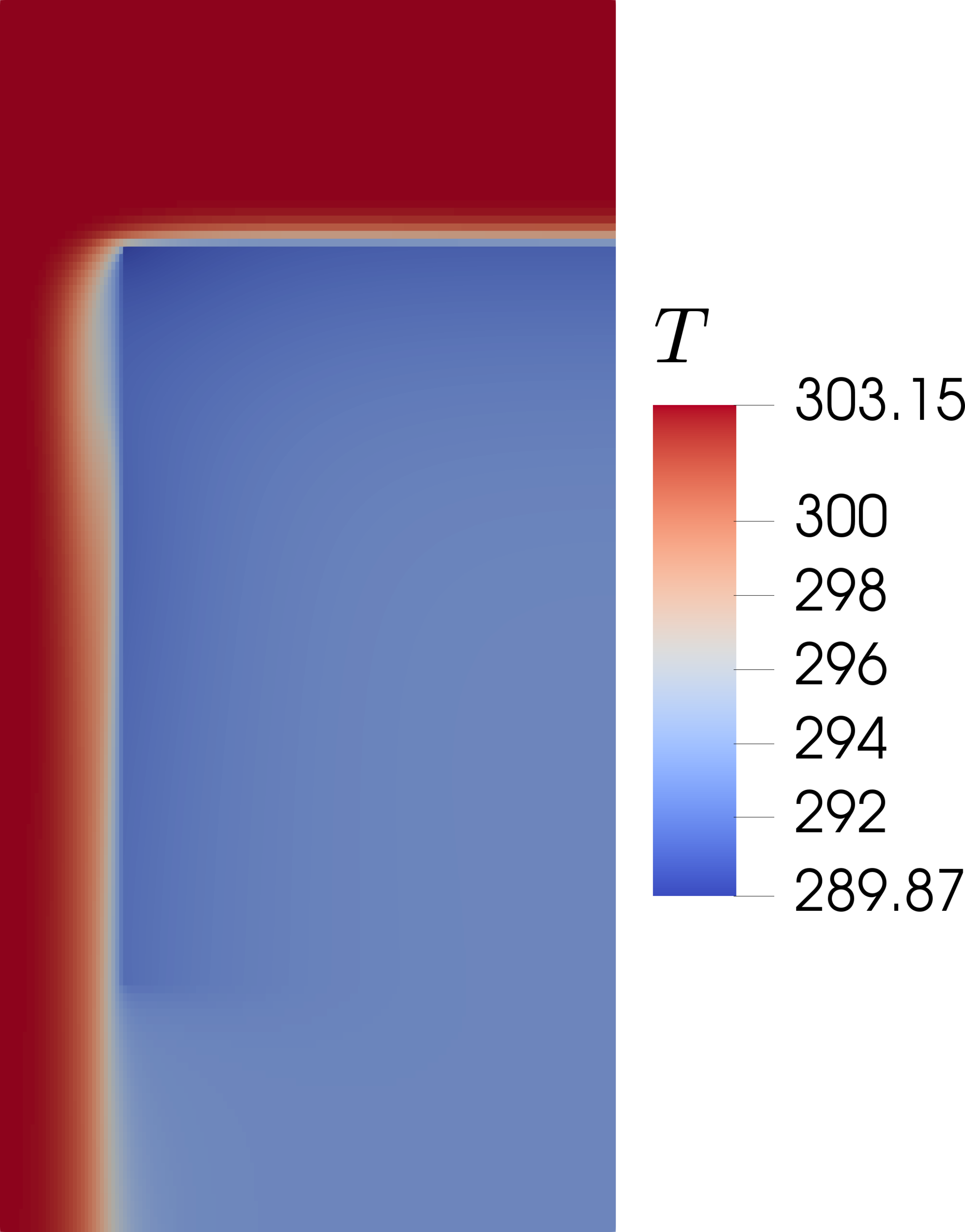}
        \caption{}
        \label{fig:fetzer_temperature}
    \end{subfigure}
    \begin{subfigure}{0.33\textwidth}
        \centering
        \includegraphics[width=0.99\textwidth]{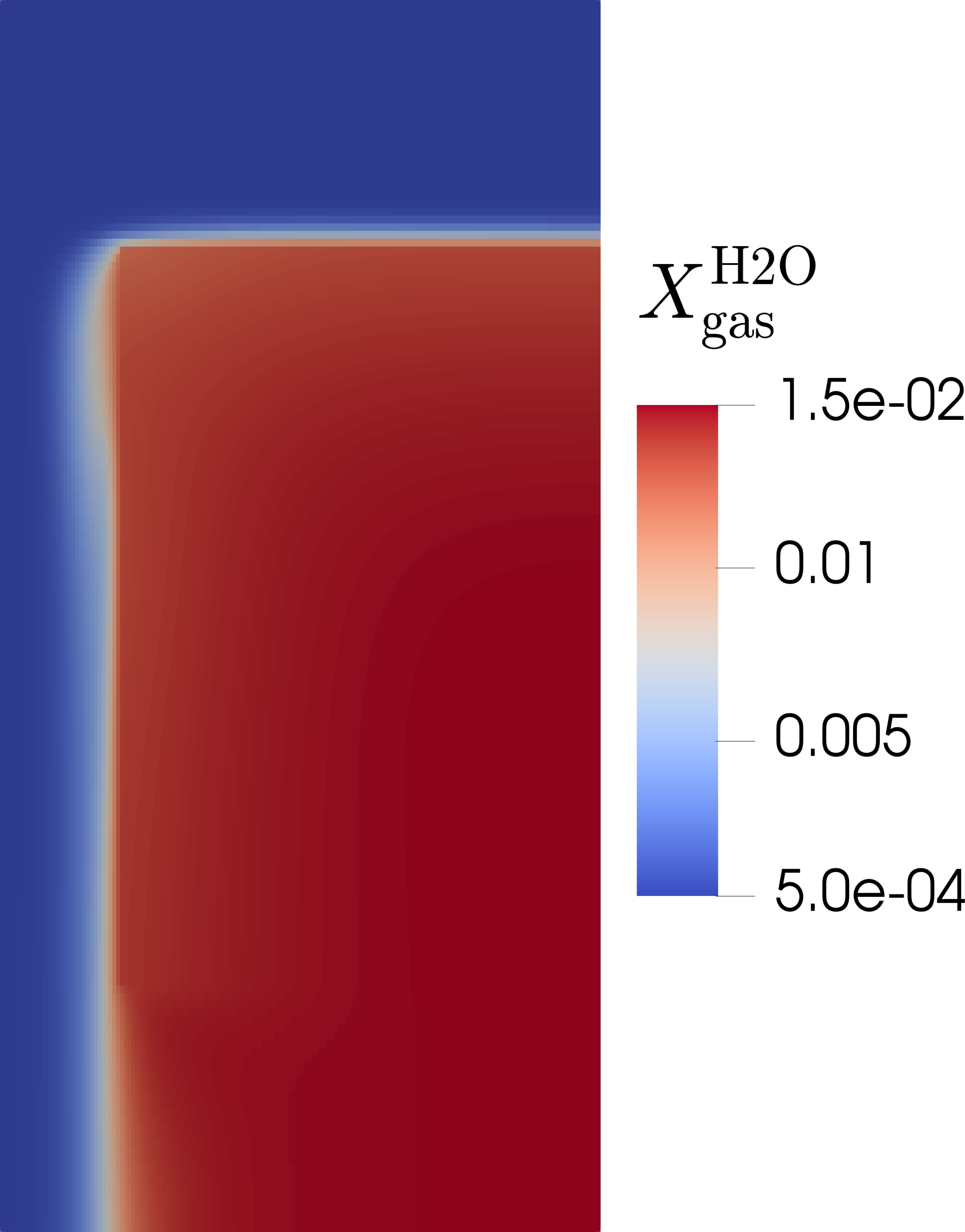}
        \caption{}
        \label{fig:fetzer_mass_frac}
    \end{subfigure}
    \begin{subfigure}{0.33\textwidth}
        \centering
        \includegraphics[width=0.99\textwidth]{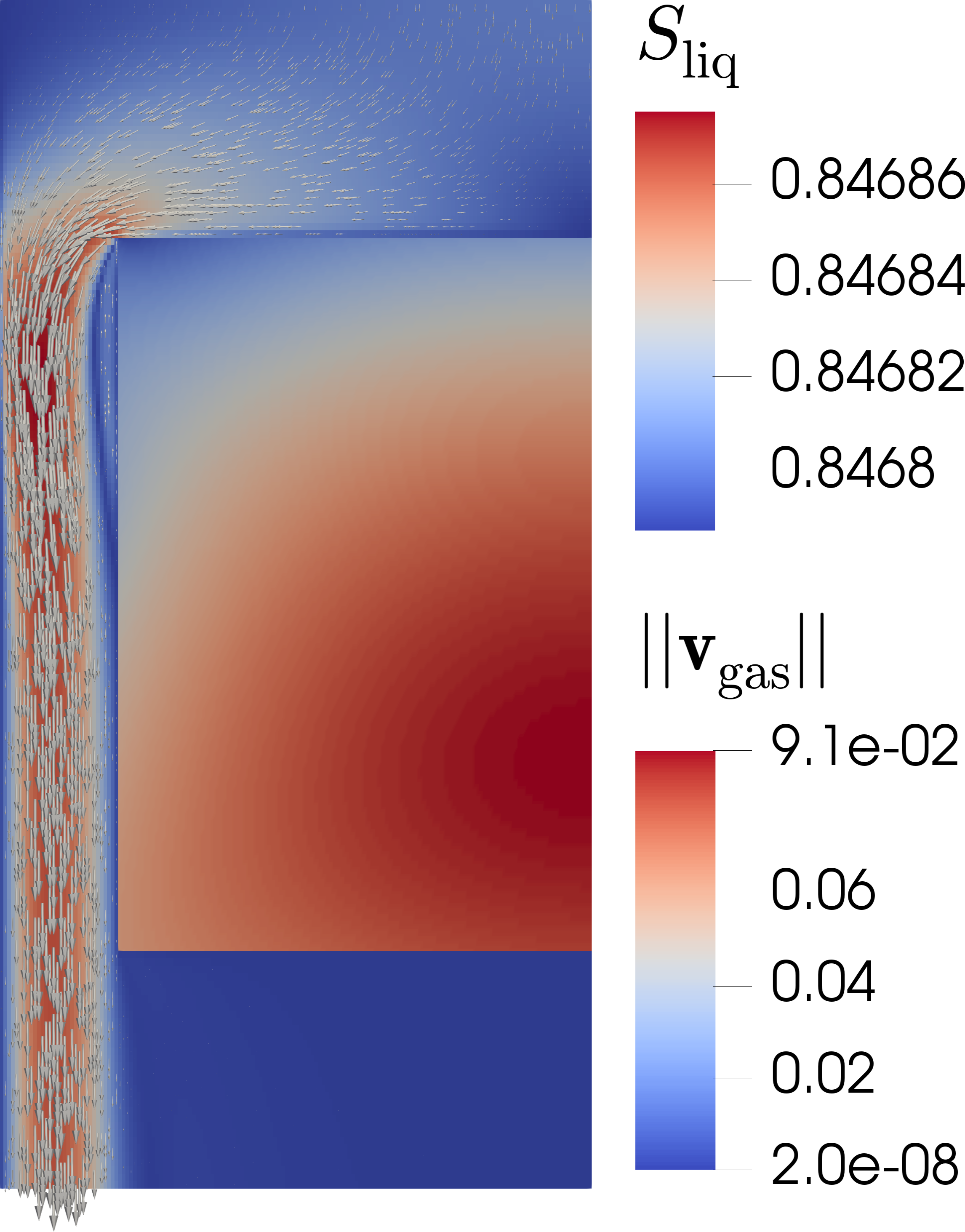}
        \caption{}
        \label{fig:fetzer_v_ff_sat_pm}
    \end{subfigure}
    \caption{\textbf{Porous obstacle}. The temperature and the mass fraction of water in the gas phase
    are depicted in (a) and (b), respectively. In (c), the velocity of air is shown in the free-flow domain while
    the water saturation is depicted in the porous medium. All snapshots were taken at the final simulation time.}
    \label{fig:fetzer_phenomenological}
\end{figure}

A two-phase, two-component model is used in the porous medium (see \cref{eq:darcy}), while a single-phase,
two-component model is considered in the free-flow domain (see \cref{eq:navierstokes}). Initially, the
porous medium is partially saturated with \emph{water} ($S_{w, \mathrm{initial}} = 0.85$) and at a
temperature of $T^\porm = 293.15 \si{\kelvin}$, while the temperature of the air in the free-flow domain
is $T^\ff = 303.15 \si{\kelvin}$. We then simulate $12 \si{\hour}$ of air flowing through the channel,
which leads to a partial desaturation of the porous medium, especially in its upper region where
the highest gas velocities occur. \Cref{fig:fetzer_phenomenological} illustrates the state at the final
simulation time. We observe that the air in the free-flow channel is cooled down by the contact and
infiltration of the cooler porous medium. Additional cooling occurs due to the evaporation of water
inside the porous medium, which leads to a decrease in temperature in the upper left region of the
porous medium, despite the fact that warmer air is flowing through it. The evaporated water is
transported away with the gas phase as depicted in the central image of \cref{fig:fetzer_phenomenological}.
The right image shows the water saturation in the porous medium, which also shows a slight decrease
towards the upper left boundary, caused by some of the water evaporating into the gas phase.

\begin{figure}
    \begin{subfigure}{0.49\textwidth}
        \centering
        \includegraphics[width=0.99\textwidth]{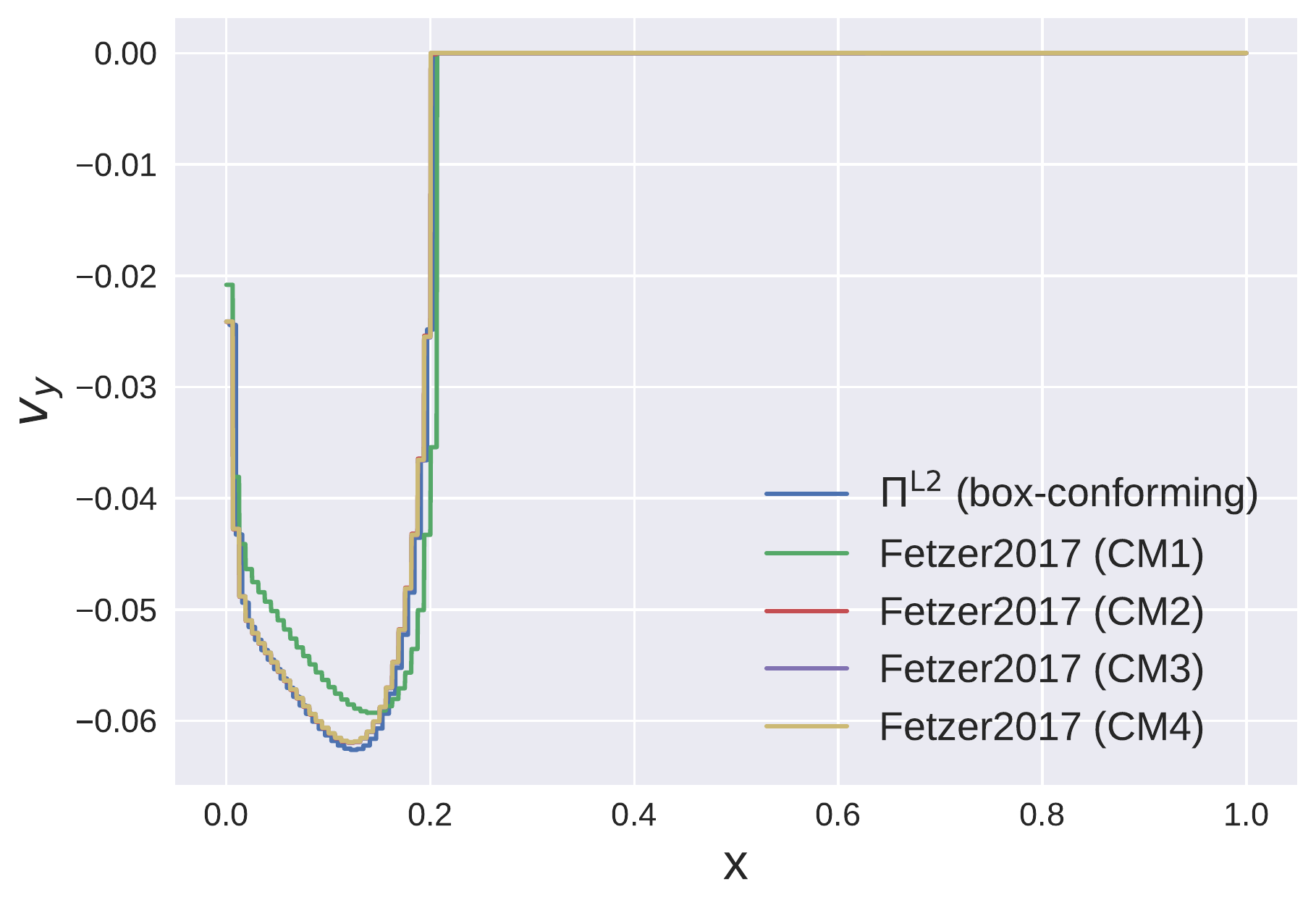}
        \caption{}
        \label{fig:fetzer_plot}
    \end{subfigure}
    \begin{subfigure}{0.49\textwidth}
        \centering
        \includegraphics[width=0.99\textwidth]{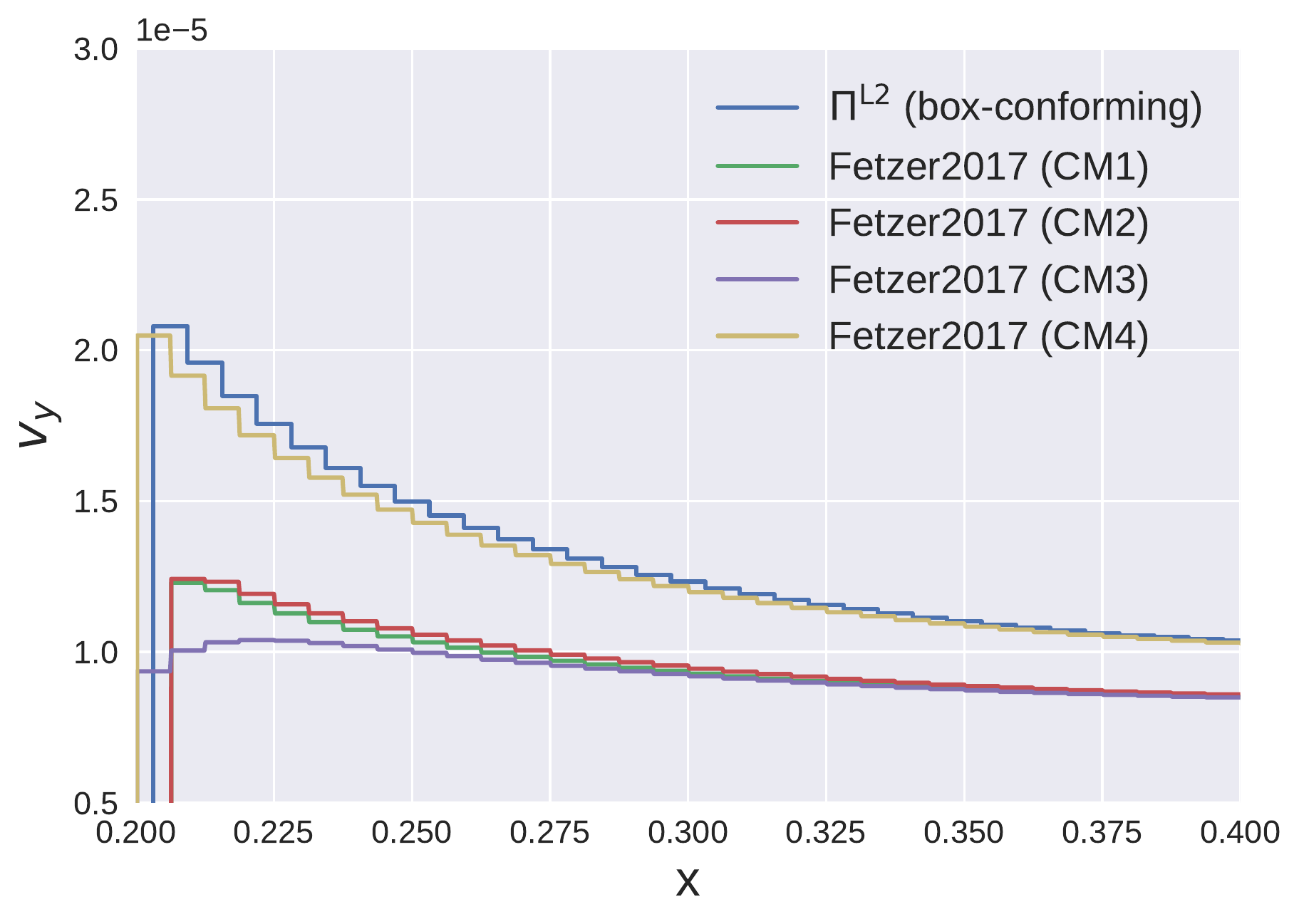}
        \caption{}
        \label{fig:fetzer_plot_zoom}
    \end{subfigure}
    \begin{subfigure}{0.49\textwidth}
        \centering
        \includegraphics[width=0.99\textwidth]{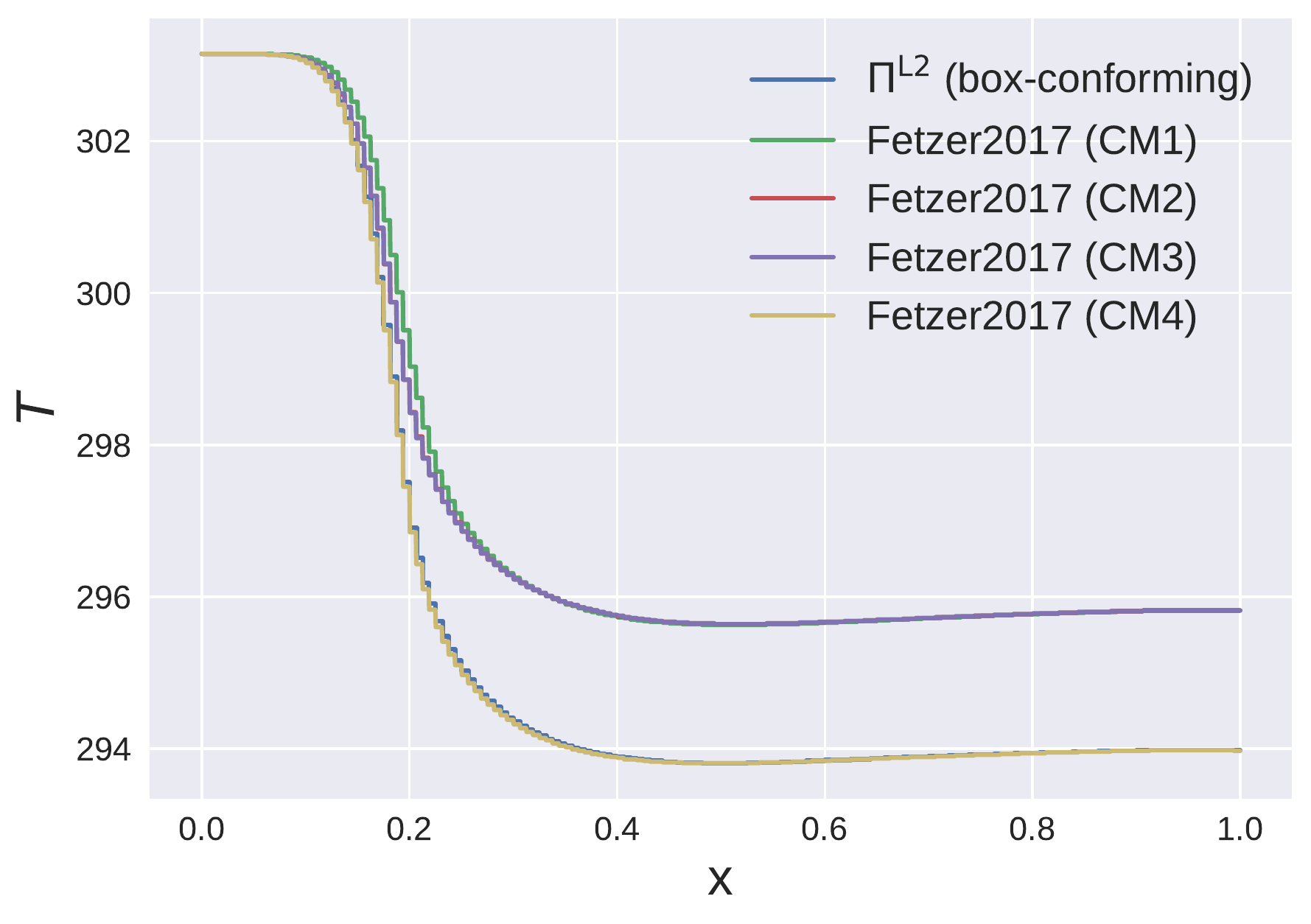}
        \caption{}
        \label{fig:fetzer_plot_temp}
    \end{subfigure}
    \begin{subfigure}{0.49\textwidth}
        \centering
        \includegraphics[width=0.99\textwidth]{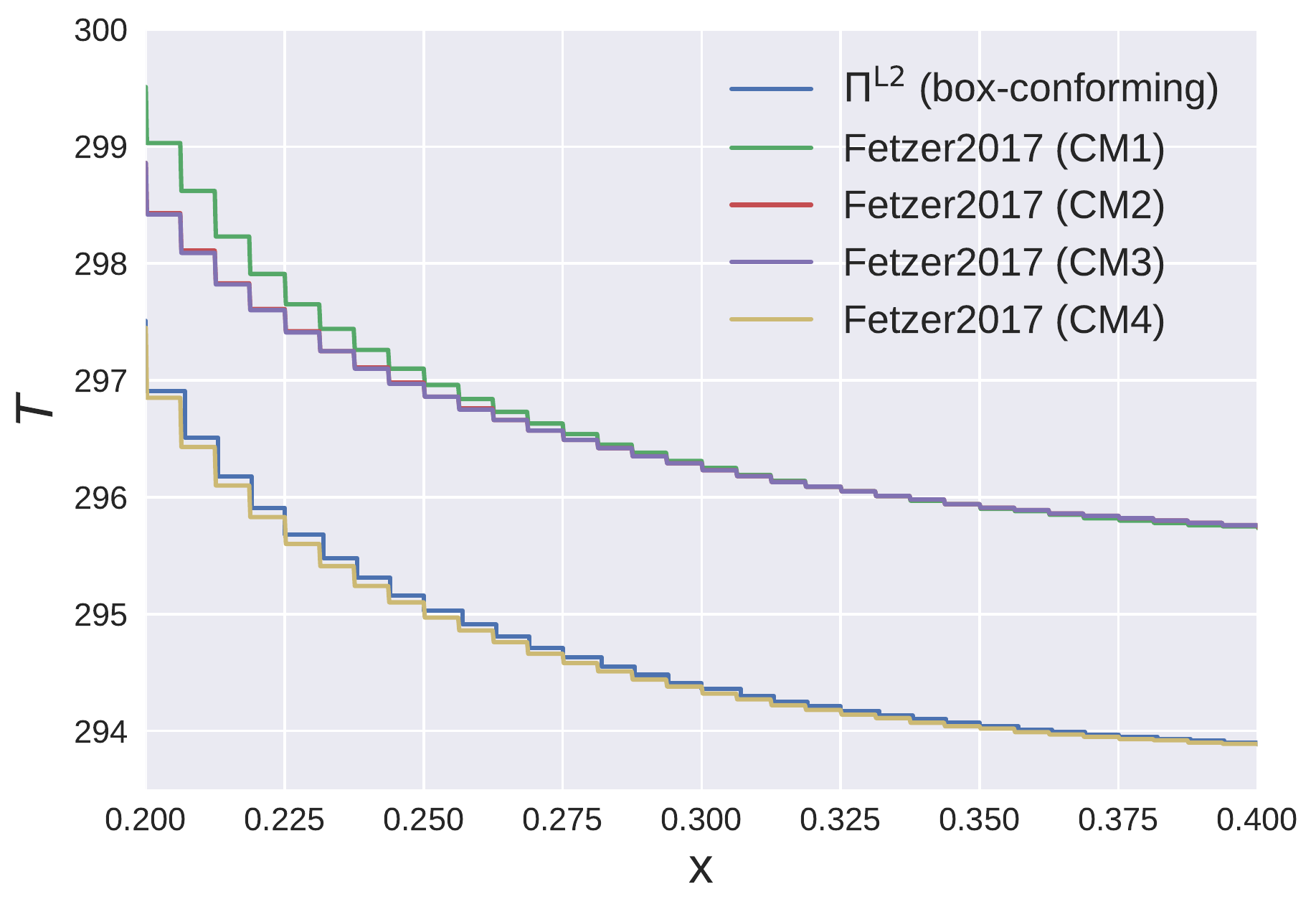}
        \caption{}
        \label{fig:fetzer_plot_temp_zoom}
    \end{subfigure}
    \caption{\textbf{Comparison to \cite{fetzer2017a}}. The figures show $v_y$ and the temperature $T$ as obtained with the proposed
    method, plotted along the line $\left(0, 1.6 \right)$ - $\left(1, 1.6 \right)$ against
    the solutions presented in \cite{fetzer2017a} with four different coupling methods. In (b) and (d), a close-up on
    the upper left part of the interface is shown.}
    \label{fig:fetzer_plots}
\end{figure}

In \cite{fetzer2017a}, four different coupling methods were compared on the basis of plots of $v_y$ along the interface.
\Cref{fig:fetzer_plots} shows the data of the original study, plotted against the results obtained with the
proposed method on a grid with the same discretization length. The coupling method ``CM4'' corresponds to the
approach that involves solving local nonlinear problems, and is considered the most accurate approach in
\cite{fetzer2017a}. The plots show that the method presented in this work leads to very similar results
without the need for an interface solver, while large deviations to the other coupling methods can be observed.

\begin{figure}
        \centering
        \includegraphics[width=0.8\textwidth]{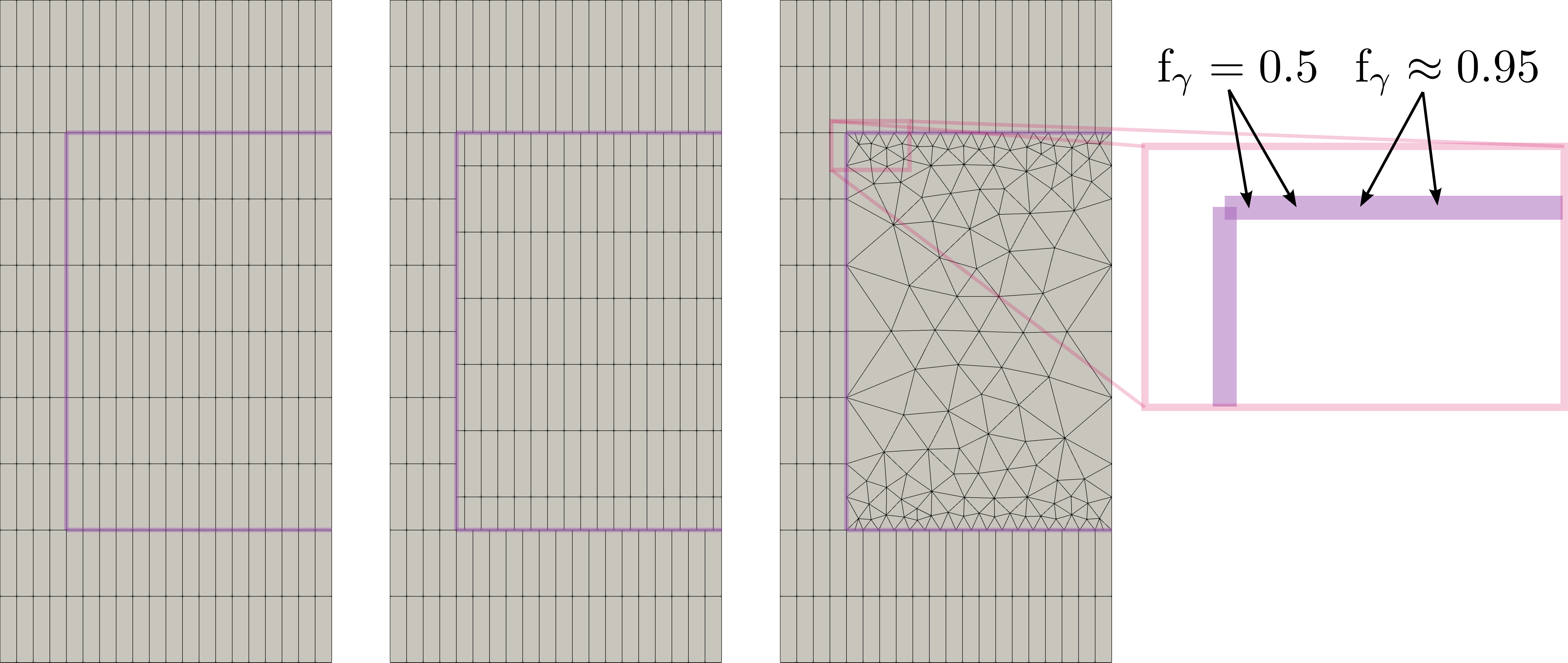}
    \caption{Conforming, box-conforming, and simplicial grid used for the test case where air flows around a porous obstacle. The interface factor $\mathrm{f}_\gamma$ is shown on the right. For the simplicial grids, faces that are located at the corners of the porous-medium domain are discretised with an interface factor of 0.5, for all other faces another factor can be set, e.g. 0.95 for the grid shown in this figure.}
    \label{fig:gridsvertical}
\end{figure}

\begin{figure}
        \centering
        \includegraphics[width=0.49\textwidth]{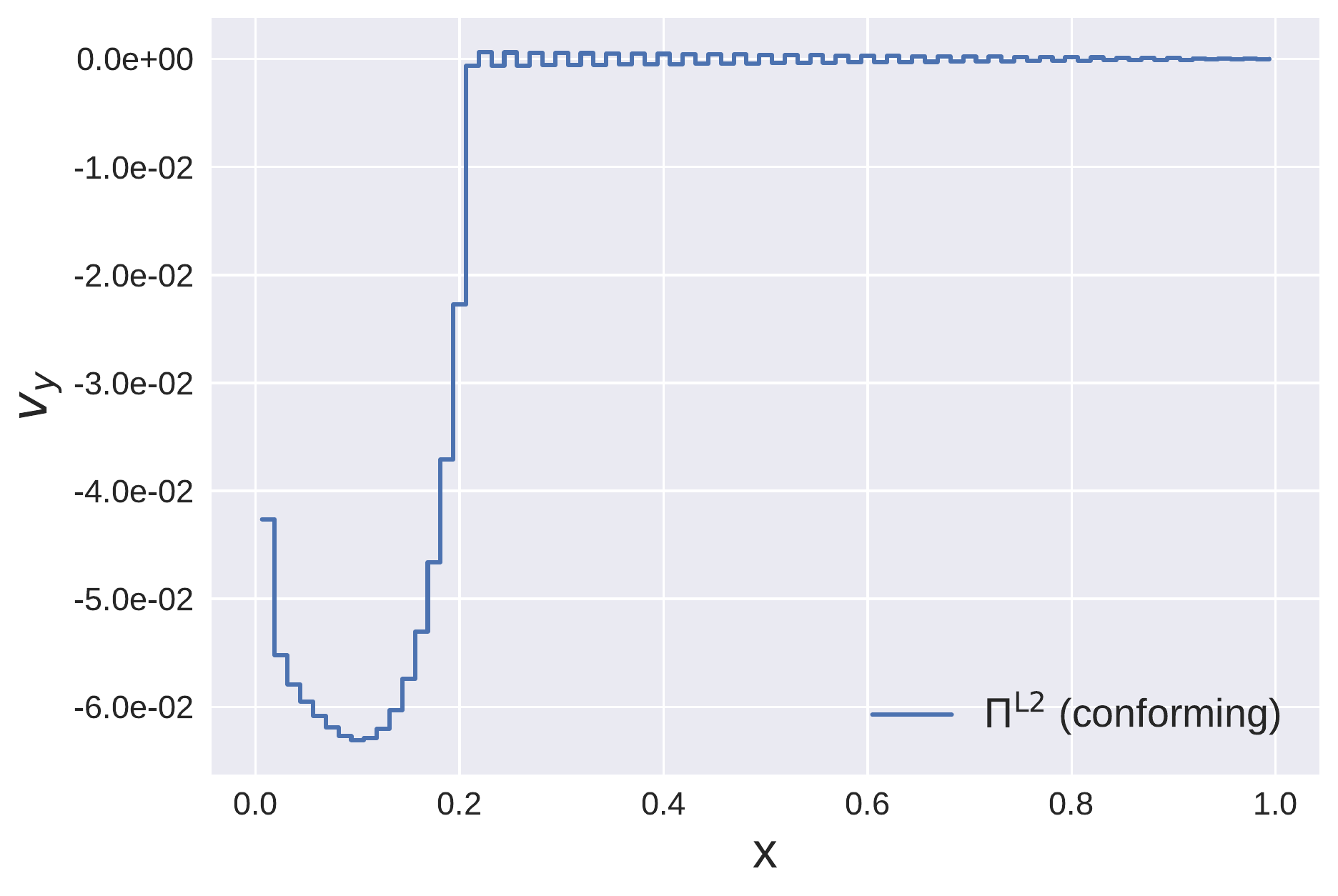}
        \includegraphics[width=0.49\textwidth]{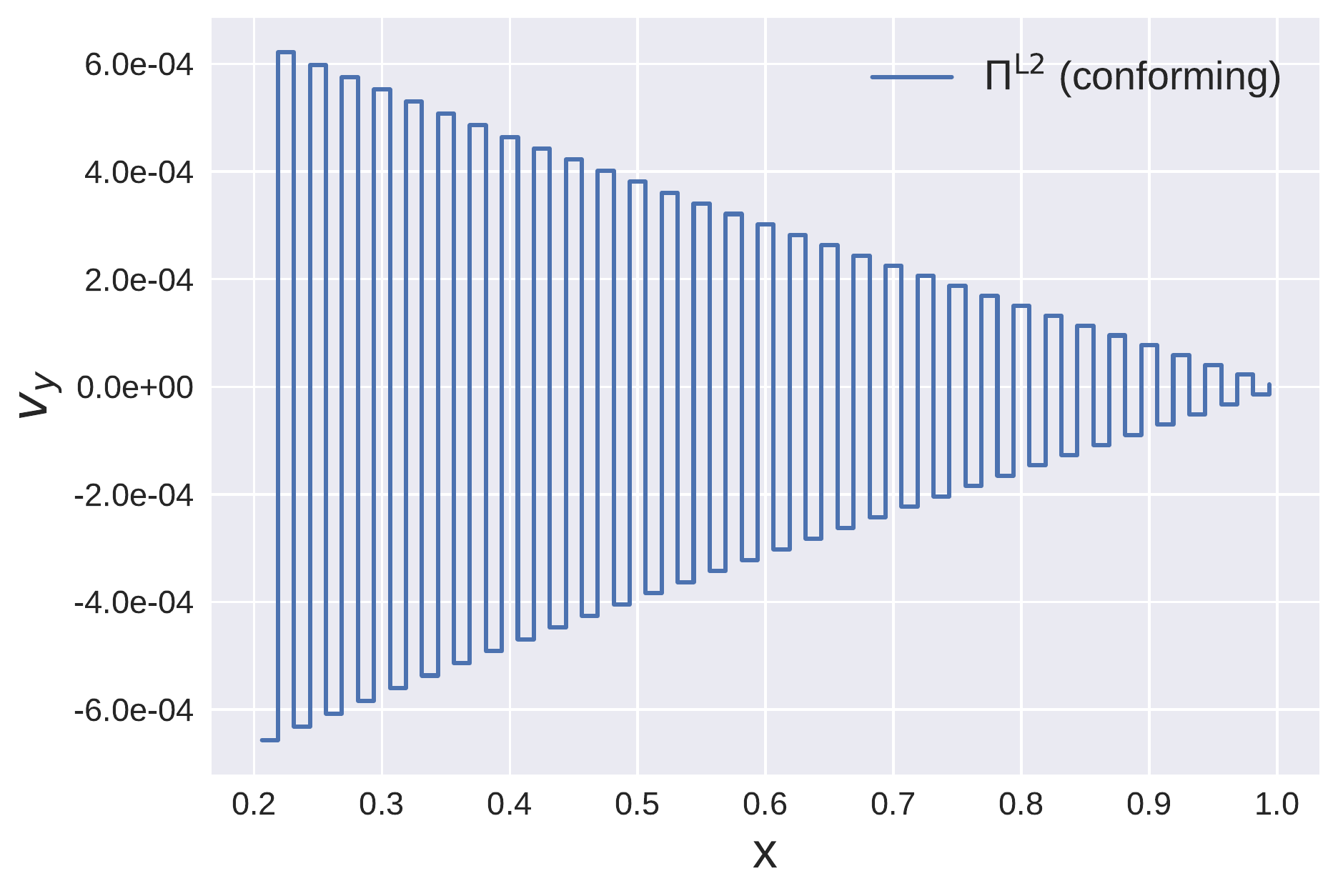}
        \includegraphics[width=0.49\textwidth]{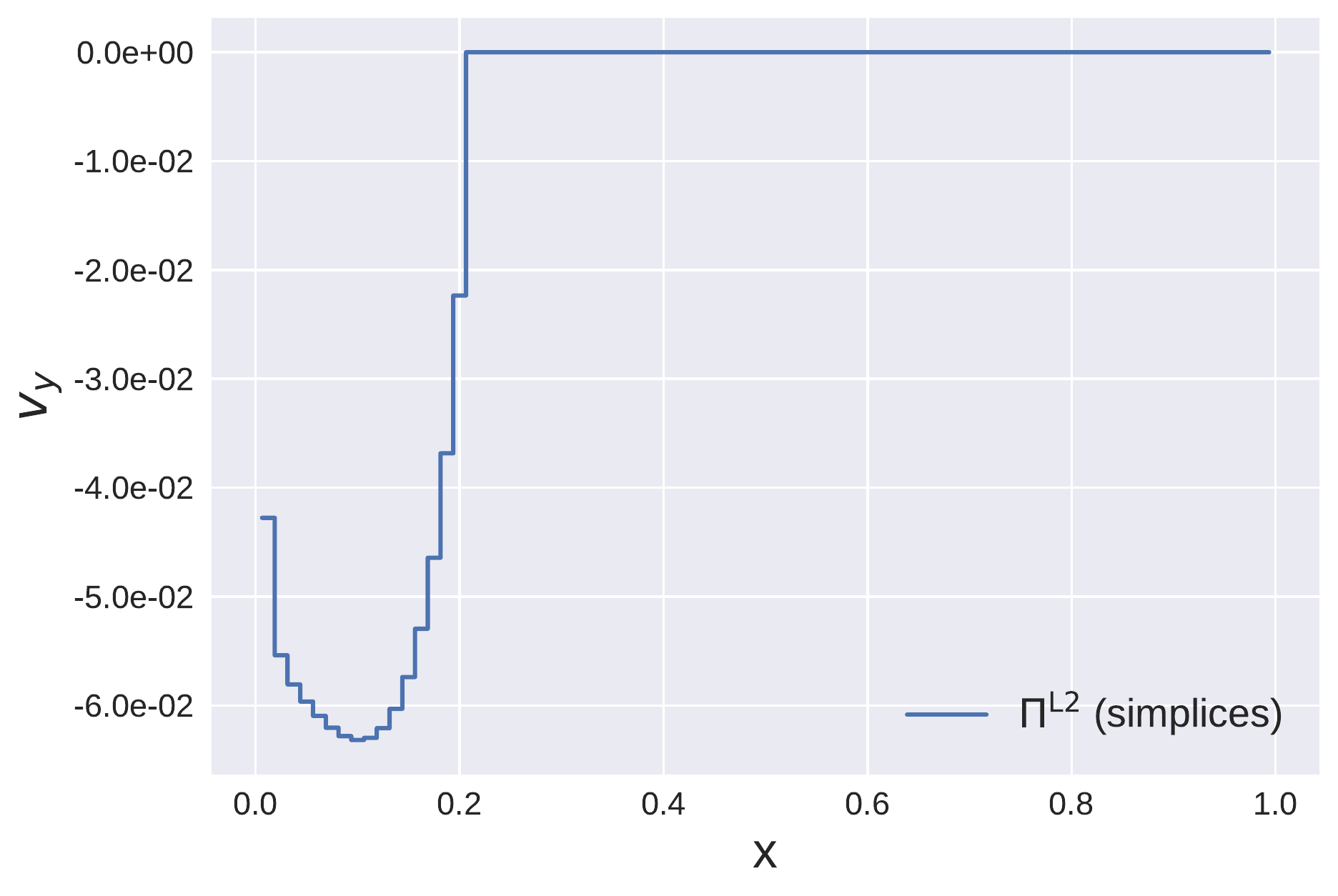}
        \includegraphics[width=0.49\textwidth]{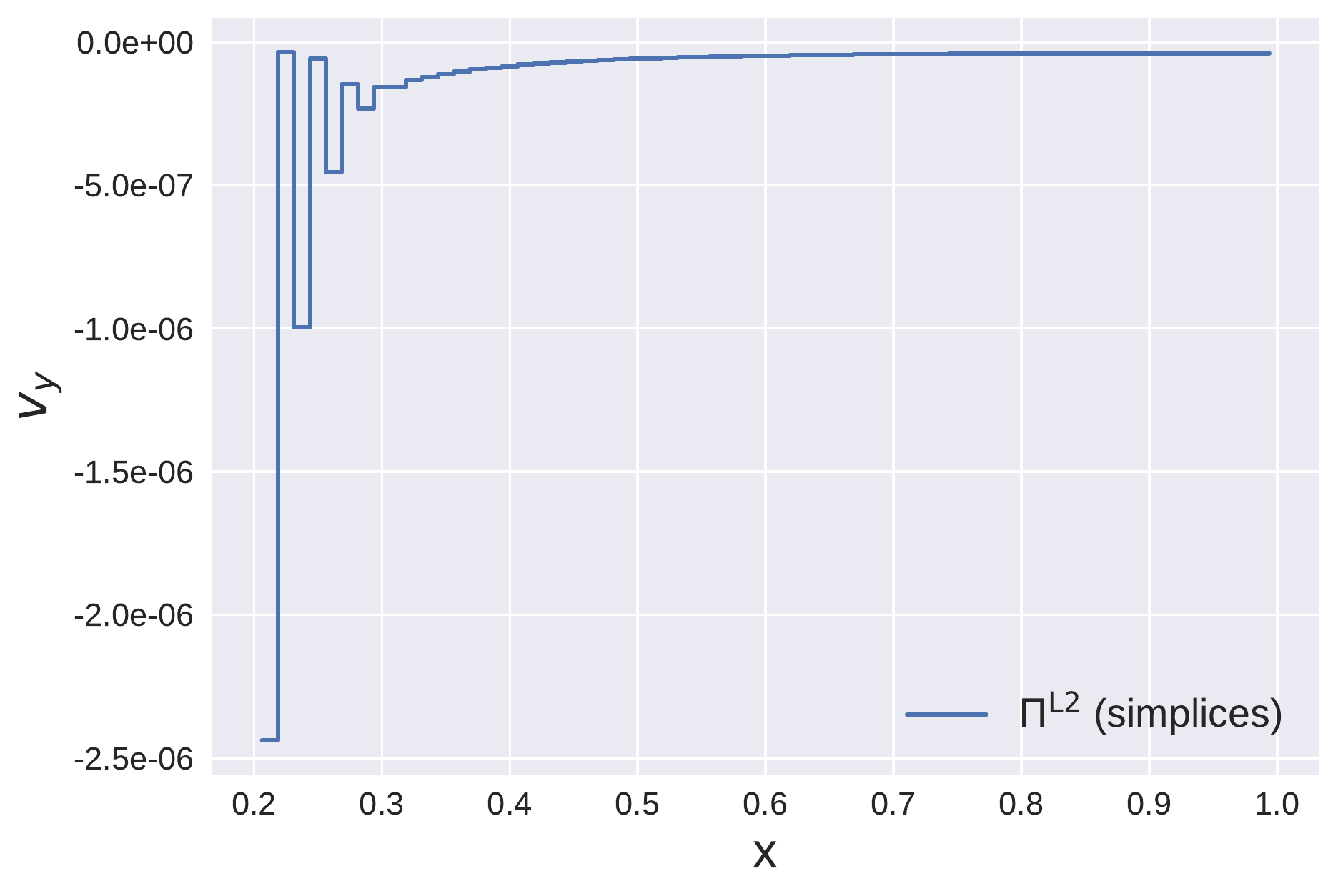}
    \caption{The figures in the left column show $v_y$ plotted along the line $\left(0, 1.6 \right)$ - $\left(1, 1.6 \right)$ obtained on a conforming grid (upper row) and on a simplicial (non-conforming) grid with an interface factor of $\mathrm{f}_\gamma \approx 0.98$.  The right column shows $v_y$ on the upper boundary of the porous medium. The discretization lengths for the free-flow grid is set to $\mathrm{d}x = \frac{1}{80}$, $\mathrm{d}y = \frac{1}{20}$ (two times refined compared to \cref{fig:gridsvertical}).}
    \label{fig:boxoscillations}
\end{figure}

\paragraph{Influence of interface grid conformity}
So far, the grid in the porous medium was chosen such that its dual grid is conforming to the grid
used in the free-flow domain. Unlike all other test cases that we have investigated so far,
we have observed an oscillatory behavior for this test case when deviating from a box-conforming grid
configuration. These oscillations also occur when neglecting non-isothermal
and compositional effects. Therefore, to avoid unnecessary complexity, we will consider a single-phase
system throughout the remainder of this section.
The corresponding results are shown in \cref{fig:boxoscillations}, where $v_y$ is plotted over the line located on the upper part of the free-flow porous-medium interface. When using a box-conforming grid (see \cref{fig:gridsvertical}) such oscillations do not occur, as discussed before and shown in \cref{fig:fetzer_plots}.

By neglecting non-isothermal and compositional effects and by considering only single-phase flow, the mass balance equation \cref{eq:darcyMassBalance} for a stationary state simplifies to $ \nabla \cdot \mathbf{\varrho v} = 0$ and
the flux coupling condition \cref{eq:ifMassFracCont} reduces to
$\left[\left(\varrho_g  \vel_g  \right)\cdot \n \right]^\ff + \left[ \left(\varrho_g \vel_g \right)\cdot \n \right]^\porm  = 0$.
Thus, the local discrete balance equation (see \cref{sub:box}) associated with each \BOX degree of freedom located at the interface is composed of interface fluxes, given by
$- |\sigma| \left[\left(\varrho_g  \vel_g  \right)\cdot \n \right]^\ff$, and of interior fluxes given by \cref{eq:discreteMassFluxes}.
This means that derivatives of these local discrete balance equations (local residuals) at the interface with respect to free-flow velocity degrees of freedom are in
the range of $\mathcal{O}(|\sigma|  \varrho_g)$ in contrast to derivatives with respect to \BOX pressure Dofs which additionally scale with the permeability $K$.
Therefore, entries in the Jacobian matrix corresponding to local \BOX residuals at the interface significantly vary.
It is therefore not surprising that numerical artifacts may occur if this is additionally amplified by grid effects. These points are further discussed in the following and the oscillatory behavior is investigated in more detail.

Similar as for the convergence test case, three grid configurations, i.e.
\emph{conforming}, \emph{non-conforming}, and \emph{box-conforming} (see \cref{fig:gridsvertical}),
are considered. For the non-conforming configuration, a triangular grid (simplices) is used in the porous-medium domain,
which allows to specify the discretization lengths at $\Gamma^\ipmff$, see \cref{fig:gridsvertical}.

To quantify the oscillations that occur at the interface (as exemplarily shown in \cref{fig:boxoscillations}), we introduce the total variaion
\begin{equation}
\mathrm{TV_\Gamma(v^\ff_y)} := \sum_{ \sigma,\sigma^\prime \subset \Gamma,
\overline{\sigma} \cap \overline{\sigma}^\prime \not=\emptyset} |v^\ff_{y,\sigma} - v^\ff_{y,\sigma^\prime}|
\end{equation}
which accumulates the differences of vertical face velocities between neighboring free-flow faces along
$\Gamma \subseteq \Gamma^\ipmff$.
In the following, we will use $\Gamma = \{(x,y) \in \Gamma^\ipmff \,|\, y = 1.6, x\geq 0.2 \}$,
that is, we only consider the upper part of the free-flow porous-medium interface. For simplicity, we omit the
subscript and denote this measure as $\mathrm{TV(v^\ff_y)}$.

\begin{figure}
    \begin{subfigure}{0.49\textwidth}
        \centering
        \includegraphics[width=0.99\textwidth]{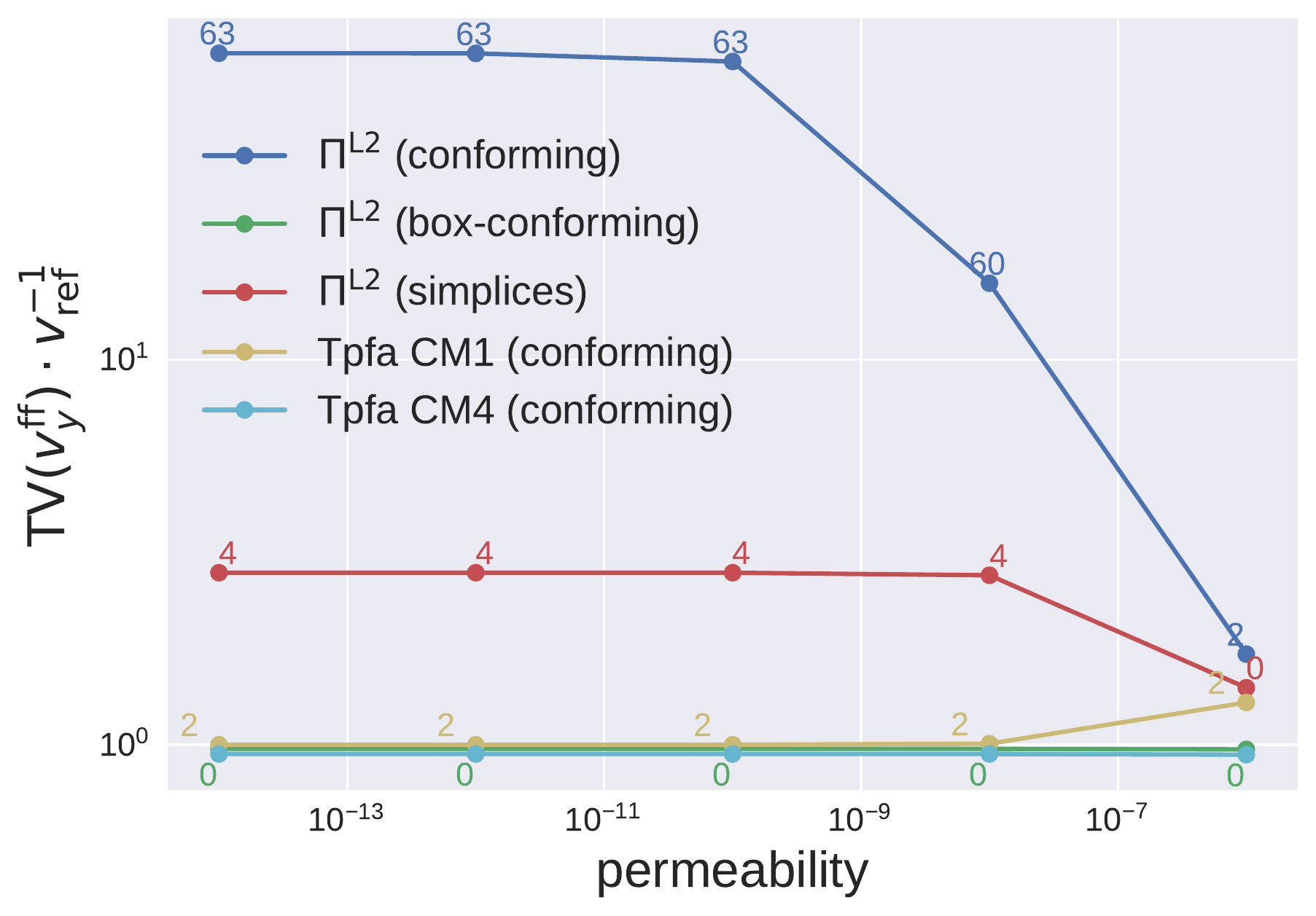}
        \caption{}
        \label{fig:convtest_grid_conforming}
    \end{subfigure}
    \begin{subfigure}{0.49\textwidth}
        \centering
        \includegraphics[width=0.99\textwidth]{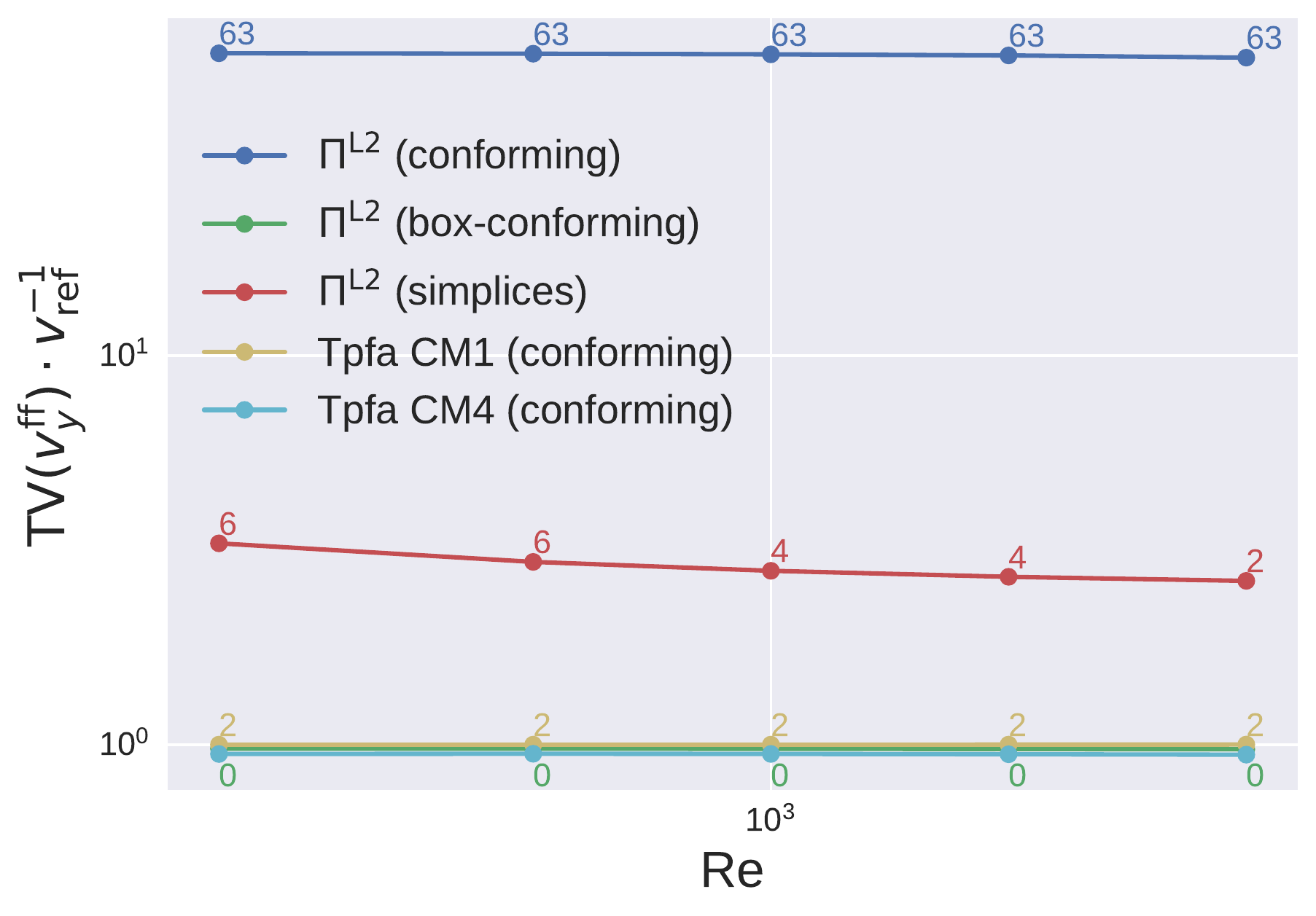}
        \caption{}
        \label{fig:convtest_grid_simplices}
    \end{subfigure}
    \caption{Total variation at upper part of the free-flow porous-medium interface, i.e. $y = 1.6$ for different permeability values and Reynolds numbers. $v_\mathrm{ref}$ is chosen as the maximum absolute velocity value that occurs on this interface. The integer values above the data points are calculated as the total number of sign changes. For the \BOX method, three grid configurations are considered, i.e. conforming, box-conforming, and non-conforming (simplices). The simplicial grid corresponds to an interface factor of $\approx 0.98$, i.e. $\sigma^\porm \approx 0.98\sigma^\ff$ for all faces $\sigma^\porm,\sigma^\ff \subset \Gamma^\ipmff$, see \cref{fig:gridsvertical}.}
    \label{fig:TVpermRe}
\end{figure}

\Cref{fig:TVpermRe} shows the total variation for different permeability values and Reynolds numbers for the grid configurations shown in \cref{fig:gridsvertical}. For comparison, the results when using the Tpfa method in the porous medium are also shown.
Here, as before, CM1 denotes the simplified interface solver, where the cell pressures are used at the interface,
whereas CM4 calculates an interface pressure from the coupling conditions by solving local problems, see \cite{fetzer2017a} for more details.
Since the Tpfa implementation currently available in \Dumux is restricted to conforming meshes, only these results are shown.
The integer values above the data points show the total number of detected sign changes,
which can be calculated as $\mathrm{TV}(0.5\mathrm{sgn}(v^\ff_y))$, where $\mathrm{sgn}$ denotes the sign function.

\begin{figure}
        \centering
        \includegraphics[width=0.49\textwidth]{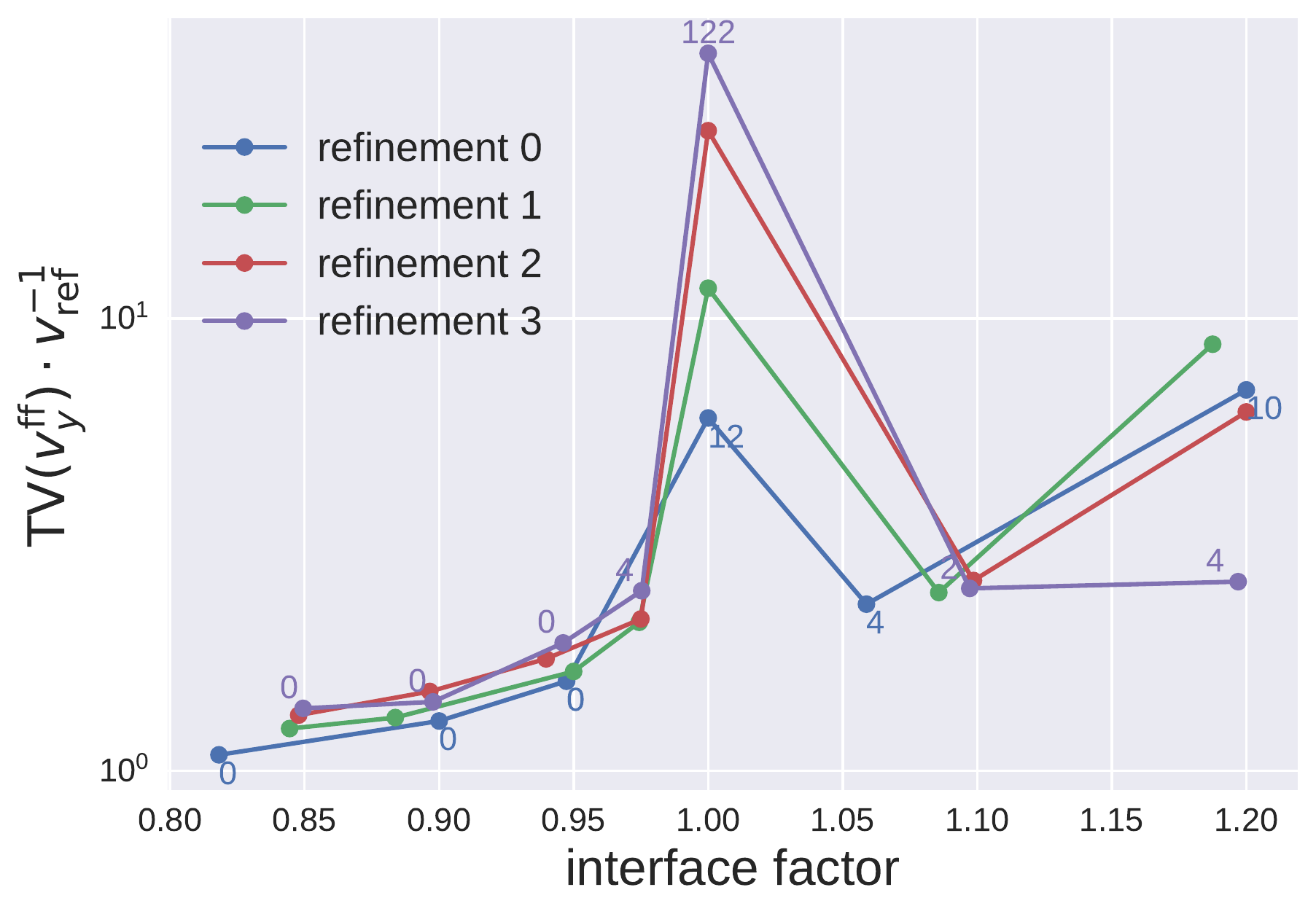}
    \caption{Total variation at upper part of the free-flow porous-medium interface, i.e. $y = 1.6$, plotted over various interface factors for different refinement levels when using the simplicial grids. $v_\mathrm{ref}$ is chosen as the maximum absolute velocity value at this interface. The integer values above the data points show the total number of detected sign changes. An interface factor of, e.g., 0.9 means that the length of a porous medium face at the interface if 0.9 times the length of a free-flow face.}
    \label{fig:dxovershoots}
\end{figure}


\Cref{fig:TVpermRe} shows that the largest total variation occurs for the conforming case with a maximum of 63 sign changes. Only for the \BOX scheme on the box-conforming grid and for the Tpfa scheme using CM4 there are no sign changes. For the Tpfa CM1 scheme, two sign changes occur close to the upper left corner of the porous medium ($(x,y) = (0.2,1.6)$). The results of the \BOX scheme for simplices with interface factor $\mathrm{f}_\gamma \approx 0.98$ is slightly worse. It can also be seen that the \BOX scheme behaves better for larger permeability values. As mentioned before, derivatives with respect to \BOX pressure Dofs directly scale with the permeability. Therefore, for larger permeability values, the magnitude of derivatives with respect to $v^\ff$ are closer to derivatives with respect to $p^\porm$.

When changing the Reynolds numbers, as shown in the right plot of \cref{fig:TVpermRe}, the behavior of the different schemes does not significantly change. This is most likely again related to the fact that derivatives of the local  \BOX residuals do not scale with the Reynolds number, because boundary fluxes linearly depend on the free-flow velocities. Note that the area-weighted projection
$\avgproj$ yielded very similar results, which are therefore not shown here.

To investigate the influence of the interface grid refinement ratio, the interface factor $\mathrm{f}_\gamma$ is varied. An interface factor of $1.0$ corresponds to the conforming case for all faces except for those that are located at the corners of the porous medium where $\mathrm{f}_\gamma = 0.5$ is set (see \cref{fig:gridsvertical}).
\Cref{fig:dxovershoots} shows the results of the \BOX scheme on the simplicial grids for different interface factors and for different grid refinement levels. Refinement 0 corresponds to discretization lengths of $\mathrm{d}x = \frac{1}{20}, \mathrm{d}y = \frac{1}{5}$ in the free-flow domain, which is
halved for each refinement. The worst results are observed for $\mathrm{f}_\gamma = 1$ (conforming case). For $\mathrm{f}_\gamma \approx 0.98$, which has also been used for the results of \cref{fig:TVpermRe}, the total variation and the number of sign changes already strongly decreases. For a factor of $ \mathrm{f}_\gamma \approx 0.95$ no sign changes are observed independent of the grid refinement level.

To summarize, by using a box-conforming grid or slightly refined grid ($\mathrm{f}_\gamma \leq 0.95$) the oscillatory behavior can be
prevented. With the \BOX method used in the porous-medium subdomain, this can be easily achieved by using unstructured grids with local refinement towards the interface.
Furthermore, it should be highlighted again that for all other test cases that we have considered so far, we did not observe any oscillatory behavior. In future work, we will investigate this in more detail and also compare our results with mortar approaches, where it is well-known that oscillations may occur if the discretization length of the mortar space is too fine compared to the coupled sub-domains.

\subsection{Example: Cooling Filter}
\label{subsec:filtertest}
In order to demonstrate the flexibility and function of the coupling system discussed above, an additional case is developed and shown here. This case highlights the model's flexibility for three dimensional geometries with non-conforming interfaces, the extension to non-isothermal models, and the simple incorporation of velocity dependent balance equations in the porous medium, here demonstrated with the Forchheimer term.

\begin{figure}[ht!]
	\centering
	\includegraphics[width=1.0\linewidth,keepaspectratio]{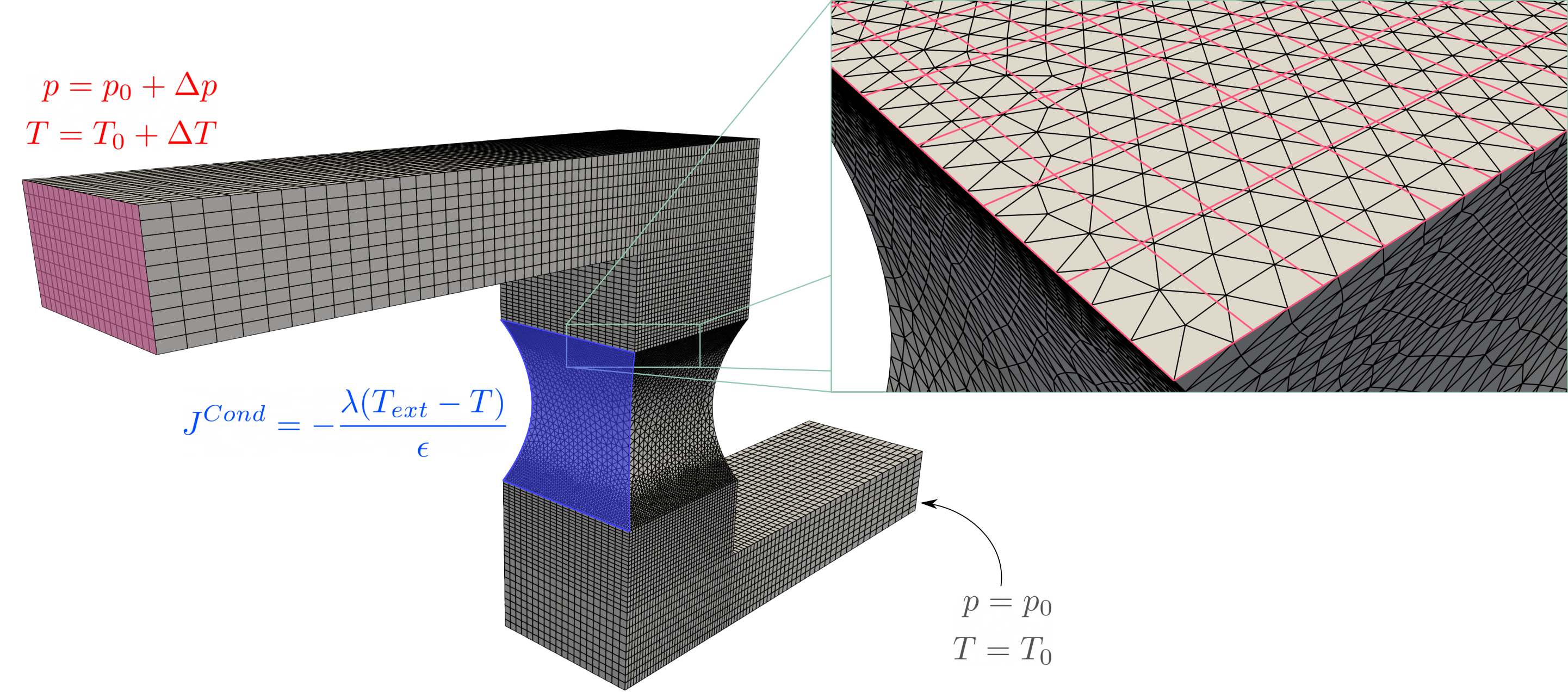}
	\caption{To the left the full simulated domain is displayed. Highlighted to the right is the difference in grid elements between the free-flow channel and the porous-medium flow section, as well as the thermal boundary conditions set at the inlet and at the left surface of the porous-medium subdomain.}
	\label{fig:filtersetup}
\end{figure}

The setup, as seen in \cref{fig:filtersetup}, incorporates two L-shaped free-flow channels, joined by a centrally located curved porous medium. Within the free-flow channels a structured grid of rectangular elements is used, where within the porous medium an unstructured tetrahedral grid is used.

A pressure difference, $\Delta P$, of $\SI{10}{\pascal}$ is applied across the system using Dirichlet boundaries for pressure at both the inlet and outlet boundaries; flow will travel through the upper L-shaped channel, through the curved porous medium, and exit through the lower L-shaped channel. All other non-interface walls within the free flow are constrained with zero mass flux and no-slip conditions. Within the porous medium, all non-interface walls are also set to allow zero mass flux. The fluid parameters simulated here are based on \Dumux's gaseous Air class defined within the H2OAir fluid system. The parameters porosity and permeability are set to $0.40$ and $\SI{1e-7}{\meter^2}$, respectively.

In addition, thermal boundary conditions are designated in the model. An increased temperature, $\Delta T$, of $\SI{20}{\kelvin}$ is set at the free-flow domain inflow, shown in red in \cref{fig:filtersetup}, as opposed to the ambient temperature, $T_0$, $\SI{298.15}{\kelvin}$. One external porous-medium boundary surface, shown in blue, will act as a heat sink. This flux is defined as a conductive flux, $J^{Cond}$ through a casing with a prescribed thickness, $\epsilon$, and conductivity $\lambda$, and with a fixed external cooling temperature $T_{ext}$. In this case, the external cooling temperature is set to be  $\SI{268.15}{\kelvin}$, the insulation thickness to $\SI{0.01}{\meter}$, and the conductivity $\SI{500}{\watt \per (\meter \kelvin)}$. All other external boundaries are set to no-temperature flux.

In addition, the Forchheimer nonlinear extension to Darcy's law \cite{nield2006a}, a velocity dependent term shown in \cref{eq:forchheimer_law}, has been added to the momentum balance within the porous medium. The permeability, $K$, is assumed to be a scalar value within this example. In this case, the coefficient $c_F$ is chosen as a fixed $0.55$ \cite{nield2006a}. Using the \BOX method, incorporation of this term is relatively simple and does not require any stencil extensions or velocity reconstructions, as would be the case for standard cell-centered finite volume methods. In this case, velocity values at the integration points are directly available from the shape functions.
\begin{equation}
 \vel_\alpha + c_F \sqrt{K} \frac{\rho_\alpha}{\mu_\alpha }
\left| \vel_\alpha\right| \vel_\alpha
+ \frac{k_{r \alpha}}{\mu_\alpha} K \nabla \left(p_\alpha
+ \rho_\alpha \mathbf{g} \right) =  0 \label{eq:forchheimer_law}
\end{equation}

\begin{figure}[ht!]
	\centering
	\includegraphics[width=1.0\linewidth,keepaspectratio]{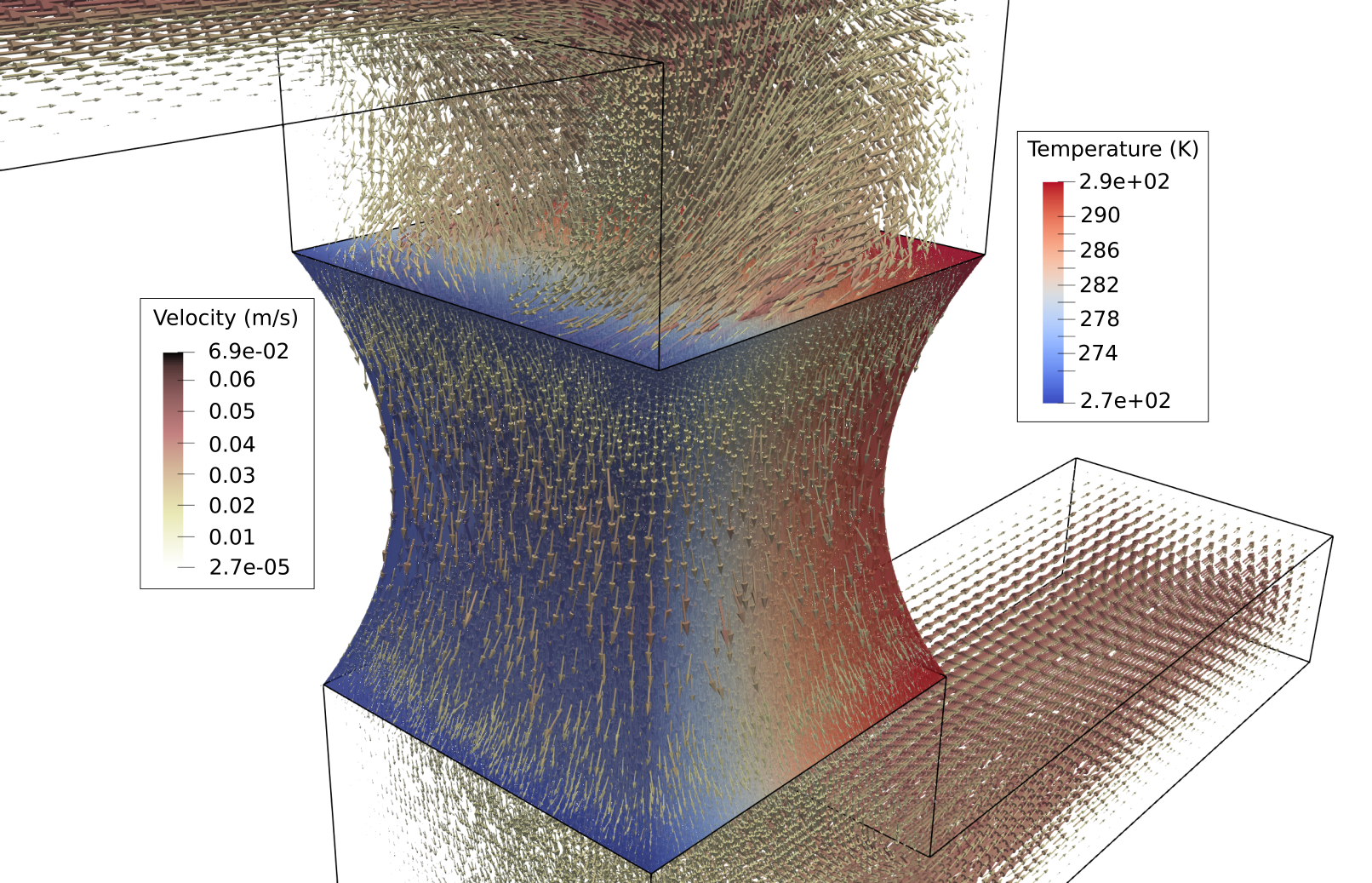}
	\caption{Flow fields are depicted in both domains using arrows, and the temperature field within the porous-medium domain is shown in the background. The cooling effect of the left cooling surface is seen to reach into the porous medium domain.}
	\label{fig:filterTempFlow}
\end{figure}

As is seen in Figure \ref{fig:filtersetup}, the interfaces between the free-flow and the porous-medium flow subdomains are non-conforming. This enables refinement to the interface and the use of unstructured grids in the porous-medium flow domain, which is required \eg to model curved domains. Seen in red is an overlay of the free-flow domain's cell boundaries at the interface. In black are the porous-medium domain's element outlines.
Across the two L-shaped free-flow domains the temperature decreases, and heat is further removed from the system across the cooling interface, as seen in Figure \ref{fig:filterTempFlow}. Although this does not represent a concrete application case, these physics are representative of applications such as CPU cooling, where heat is removed from a porous system embedded within a larger flow system \cite{Zhang2018a}.

\section{Conclusions}
\label{sec:conclusions}
This work explores coupled two-domain free-flow and multi-phase porous medium flow models using a vertex-centered finite volume
method to discretize the porous-medium flow domain.
Use of this \BOX discretization method in the porous-medium flow subdomain allows for the use of unstructured non-conforming grids, more flexible inclusion of interfacial coupling conditions, and a simpler incorporation of velocity dependent terms in the balance equations such as the Forchheimer extension of Darcy's law, as shown in this work.

The convergence behavior of the proposed method is investigated in \cref{subsec:convtest}. Therein, conforming, non-conforming and box-conforming grid configurations are considered and convergence orders are determined across levels of refinement in comparison to an analytical solution for two projection operators \ltwoproj and \avgproj. For both projection operators
and all considered grids, second order convergence is observed.

The performance of the proposed method in a non-isothermal, multi-component and multi-phase environment is examined in \ref{sec:fetzer}. The accuracy of the coupling system was then further investigated to provide insight into limitations due to oscillations when dealing with non box-conforming unstructured grids. Recommendations for grid constraints are made to ensure non-oscillatory behavior. In practice, it suffices to locally refine the grid used in the porous medium towards the interface.

Finally, an exemplary 3D application is presented, which illustrates the capabilities of the proposed method:
The porous-medium subdomain is curved and is discretized with an unstructured tetrahedral grid, while the interfaces
to the neighboring free-flow channels are non-conforming. Moreover, within the porous-medium subdomain, the nonlinear velocity dependent Forchheimer extension of Darcy's law is incorporated, and heat transport is evaluated across the full coupled system.

As the use and development of coupled porous medium free-flow models are a field of active research, future work in this area will include both extensions to the discretization methods as well as developments to the incorporated mathematical models.
At the domain interface, further coupling conditions should be investigated to better resolve the near interface conditions \cite{eggenweiler2021a}.
Additionally, at the moment the two domain system is solved monolithically.
Comparing the existing monolithically solved method with a partitioned model where domains are solved individually \cite{jaust2020a}, or to other domain decomposition models using interface mortars \cite{boon2020a}, could also provide insight to other solution methods.
For expansions to the free-flow domain, turbulence models and radiation based expansions to the energy balance could be introduced to better simulate atmospheric free-flow conditions \cite{heck2020a, coltman2020a}. To expand on the discretization in the free-flow domain, unstructured and adaptive staggered finite volume methods are being developed in order to efficiently resolve more complex geometries and flow cases.

\section*{Acknowledgements}
We thank the German Research Foundation (DFG) for supporting this work by funding SFB 1313, Project Number 327154368, Research Project A02,
and by funding SimTech via Germany's Excellence Strategy (EXC 2075 – 390740016).

\appendix
\section{}
\label{sec:appendix_conv_tables}

The following tables list the errors obtained in the convergence test of \Cref{subsec:convtest}
on the different grids and different projections used.

\begin{table}[h]
  \centering
  \caption{\textbf{Convergence study}. Errors and rates for the primary variables (conforming case, l2-projection).}
  \label{tab:convergence_test_conforming}
  \begin{tabular}{*{7}{l}|*2{l}}
    \toprule
     & \multicolumn{6}{c|}{free flow} &  \multicolumn{2}{c}{darcy flow} \\
     $m$ & $e_p^m$ & $r_p$ & $e_{v_x}^m$ & $r_{v_x}$ & $e_{v_y}^m$ & $r_{v_y}$ & $e_{p^\mathrm{pm}}^m$ & $r_{p^\mathrm{pm}}$ \\
    \midrule
    0 & 1.44e-01 & - & 2.54e-03 & - & 8.97e-03 & - & 2.52e-02 & - \\
    1 & 3.76e-02 & 1.94e+00 & 5.35e-04 & 2.25e+00 & 2.10e-03 & 2.10e+00 & 6.11e-03 & 2.05e+00 \\
    2 & 9.52e-03 & 1.98e+00 & 1.27e-04 & 2.07e+00 & 5.15e-04 & 2.03e+00 & 1.50e-03 & 2.03e+00 \\
    3 & 2.39e-03 & 1.99e+00 & 3.15e-05 & 2.02e+00 & 1.28e-04 & 2.01e+00 & 3.72e-04 & 2.01e+00 \\
    4 & 5.98e-04 & 2.00e+00 & 7.84e-06 & 2.00e+00 & 3.20e-05 & 2.00e+00 & 9.26e-05 & 2.00e+00 \\
    5 & 1.49e-04 & 2.00e+00 & 1.96e-06 & 2.00e+00 & 7.99e-06 & 2.00e+00 & 2.31e-05 & 2.00e+00 \\
  \bottomrule
  \end{tabular}
\end{table}

\begin{table}[h]
  \centering
  \caption{\textbf{Convergence study}. Errors and rates for the primary variables (conforming case, weighted dof-projection).}
  \label{tab:convergence_test_conforming_weighteddof}
  \begin{tabular}{*{7}{l}|*2{l}}
    \toprule
     & \multicolumn{6}{c|}{free flow} &  \multicolumn{2}{c}{darcy flow} \\
     $m$ & $e_p^m$ & $r_p$ & $e_{v_x}^m$ & $r_{v_x}$ & $e_{v_y}^m$ & $r_{v_y}$ & $e_{p^\mathrm{pm}}^m$ & $r_{p^\mathrm{pm}}$ \\
    \midrule
    0 & 1.44e-01 & - & 2.54e-03 & - & 8.97e-03 & - & 2.52e-02 & - \\
    1 & 3.76e-02 & 1.94e+00 & 5.35e-04 & 2.25e+00 & 2.10e-03 & 2.10e+00 & 6.11e-03 & 2.05e+00 \\
    2 & 9.52e-03 & 1.98e+00 & 1.27e-04 & 2.07e+00 & 5.15e-04 & 2.03e+00 & 1.50e-03 & 2.03e+00 \\
    3 & 2.39e-03 & 1.99e+00 & 3.15e-05 & 2.02e+00 & 1.28e-04 & 2.01e+00 & 3.72e-04 & 2.01e+00 \\
    4 & 5.98e-04 & 2.00e+00 & 7.84e-06 & 2.00e+00 & 3.20e-05 & 2.00e+00 & 9.26e-05 & 2.00e+00 \\
    5 & 1.49e-04 & 2.00e+00 & 1.96e-06 & 2.00e+00 & 7.99e-06 & 2.00e+00 & 2.31e-05 & 2.00e+00 \\
  \bottomrule
  \end{tabular}
\end{table}

\begin{table}[h]
  \centering
  \caption{\textbf{Convergence study}. Errors and rates for the primary variables (simplices, l2-projection).}
  \label{tab:convergence_test_simplices}
  \begin{tabular}{*{7}{l}|*2{l}}
    \toprule
     & \multicolumn{6}{c|}{free flow} &  \multicolumn{2}{c}{darcy flow} \\
     $m$ & $e_p^m$ & $r_p$ & $e_{v_x}^m$ & $r_{v_x}$ & $e_{v_y}^m$ & $r_{v_y}$ & $e_{p^\mathrm{pm}}^m$ & $r_{p^\mathrm{pm}}$ \\
    \midrule
    0 & 1.64e-01 & - & 2.29e-03 & - & 8.69e-03 & - & 2.61e-02 & - \\
    1 & 4.64e-02 & 1.82e+00 & 4.67e-04 & 2.29e+00 & 2.02e-03 & 2.10e+00 & 6.58e-03 & 1.99e+00 \\
    2 & 1.21e-02 & 1.94e+00 & 1.08e-04 & 2.11e+00 & 4.89e-04 & 2.05e+00 & 1.66e-03 & 1.99e+00 \\
    3 & 3.11e-03 & 1.96e+00 & 2.69e-05 & 2.01e+00 & 1.21e-04 & 2.01e+00 & 4.17e-04 & 1.99e+00 \\
    4 & 7.81e-04 & 2.00e+00 & 6.67e-06 & 2.01e+00 & 3.02e-05 & 2.01e+00 & 1.04e-04 & 2.00e+00 \\
    5 & 1.96e-04 & 1.99e+00 & 1.68e-06 & 1.99e+00 & 7.56e-06 & 2.00e+00 & 2.61e-05 & 2.00e+00 \\
  \bottomrule
  \end{tabular}
\end{table}

\begin{table}[h]
  \centering
  \caption{\textbf{Convergence study}. Errors and rates for the primary variables (simplices, weighted dof-projection).}
  \label{tab:convergence_test_simplices_weighteddof}
  \begin{tabular}{*{7}{l}|*2{l}}
    \toprule
     & \multicolumn{6}{c|}{free flow} &  \multicolumn{2}{c}{darcy flow} \\
     $m$ & $e_p^m$ & $r_p$ & $e_{v_x}^m$ & $r_{v_x}$ & $e_{v_y}^m$ & $r_{v_y}$ & $e_{p^\mathrm{pm}}^m$ & $r_{p^\mathrm{pm}}$ \\
    \midrule
    0 & 1.53e-01 & - & 2.19e-03 & - & 8.26e-03 & - & 2.61e-02 & - \\
    1 & 4.35e-02 & 1.81e+00 & 4.60e-04 & 2.25e+00 & 2.01e-03 & 2.04e+00 & 6.58e-03 & 1.99e+00 \\
    2 & 1.12e-02 & 1.96e+00 & 1.06e-04 & 2.12e+00 & 4.85e-04 & 2.05e+00 & 1.66e-03 & 1.99e+00 \\
    3 & 2.92e-03 & 1.94e+00 & 2.62e-05 & 2.01e+00 & 1.20e-04 & 2.01e+00 & 4.17e-04 & 1.99e+00 \\
    4 & 7.29e-04 & 2.00e+00 & 6.35e-06 & 2.04e+00 & 2.97e-05 & 2.02e+00 & 1.04e-04 & 2.00e+00 \\
    5 & 1.84e-04 & 1.98e+00 & 1.62e-06 & 1.97e+00 & 7.49e-06 & 1.99e+00 & 2.62e-05 & 2.00e+00 \\
  \bottomrule
  \end{tabular}
\end{table}

\begin{table}[h]
  \centering
  \caption{\textbf{Convergence study}. Errors and rates for the primary variables (box-conforming case, 2-projection).}
  \label{tab:convergence_test_box_conforming}
  \begin{tabular}{*{7}{l}|*2{l}}
    \toprule
     & \multicolumn{6}{c|}{free flow} &  \multicolumn{2}{c}{darcy flow} \\
     $m$ & $e_p^m$ & $r_p$ & $e_{v_x}^m$ & $r_{v_x}$ & $e_{v_y}^m$ & $r_{v_y}$ & $e_{p^\mathrm{pm}}^m$ & $r_{p^\mathrm{pm}}$ \\
    \midrule
    0 & 1.30e-01 & - & 2.43e-03 & - & 8.80e-03 & - & 2.16e-02 & - \\
    1 & 3.22e-02 & 2.01e+00 & 5.17e-04 & 2.24e+00 & 2.07e-03 & 2.09e+00 & 5.74e-03 & 1.91e+00 \\
    2 & 7.88e-03 & 2.03e+00 & 1.26e-04 & 2.04e+00 & 5.11e-04 & 2.02e+00 & 1.48e-03 & 1.96e+00 \\
    3 & 1.94e-03 & 2.02e+00 & 3.17e-05 & 1.99e+00 & 1.28e-04 & 2.00e+00 & 3.74e-04 & 1.98e+00 \\
    4 & 4.81e-04 & 2.01e+00 & 7.99e-06 & 1.99e+00 & 3.20e-05 & 2.00e+00 & 9.44e-05 & 1.99e+00 \\
    5 & 1.20e-04 & 2.01e+00 & 2.01e-06 & 1.99e+00 & 8.01e-06 & 2.00e+00 & 2.37e-05 & 1.99e+00 \\
  \bottomrule
  \end{tabular}
\end{table}

\begin{table}[h]
  \centering
  \caption{\textbf{Convergence study}. Errors and rates for the primary variables (box-conforming case, weighted dof-projection).}
  \label{tab:convergence_test_box_conforming_weighteddof}
  \begin{tabular}{*{7}{l}|*2{l}}
    \toprule
     $m$ & $e_p^m$ & $r_p$ & $e_{v_x}^m$ & $r_{v_x}$ & $e_{v_y}^m$ & $r_{v_y}$ & $e_{p^\mathrm{pm}}^m$ & $r_{p^\mathrm{pm}}$ \\
    \midrule
    0 & 1.04e-01 & - & 2.66e-03 & - & 8.17e-03 & - & 2.14e-02 & - \\
    1 & 2.24e-02 & 2.21e+00 & 5.79e-04 & 2.20e+00 & 1.97e-03 & 2.05e+00 & 5.71e-03 & 1.90e+00 \\
    2 & 4.96e-03 & 2.17e+00 & 1.36e-04 & 2.09e+00 & 4.95e-04 & 2.00e+00 & 1.47e-03 & 1.96e+00 \\
    3 & 1.14e-03 & 2.12e+00 & 3.29e-05 & 2.05e+00 & 1.25e-04 & 1.99e+00 & 3.73e-04 & 1.98e+00 \\
    4 & 2.71e-04 & 2.08e+00 & 8.10e-06 & 2.02e+00 & 3.13e-05 & 1.99e+00 & 9.42e-05 & 1.99e+00 \\
    5 & 6.53e-05 & 2.05e+00 & 2.02e-06 & 2.01e+00 & 7.86e-06 & 1.99e+00 & 2.37e-05 & 1.99e+00 \\
  \bottomrule
  \end{tabular}
\end{table}

\clearpage
\bibliography{literature}

\end{document}